\documentclass[aop,MSNbibl,seceqn,dvips]{arximspdf}

\doi{10.1214/12-AOP775} 
\volume{41}
\issue{2}
\pubyear{2013}
\firstpage{1088}
\lastpage{1114}

\makeatletter
\newcommand{\rrVert}{\Vert}
\newcommand{\rrvert}{\vert}
\newcommand{\llVert}{\Vert}
\newcommand{\llvert}{\vert}
\newtheorem{theorem}{Theorem}[section]
\newproclaim{definition}{Definition}[section]
\newtheorem{proposition}{Proposition}[section]
\newtheorem{corollary}{Corollary}[section]
\newtheorem{lemma}{Lemma}[section]
\newproclaim{example}{Example}
\makeatother

\begin{document}
\begin{frontmatter}

\title{Limit theorems for weighted nonlinear transformations of Gaussian
stationary processes with singular spectra}
\runtitle{Limit theorems for nonlinear transformations}

\begin{aug}
\author[A]{\fnms{Alexander V.} \snm{Ivanov}\ead[label=e1]{alexntuu@gmail.com}},
\author[B]{\fnms{Nikolai} \snm{Leonenko}\thanksref{t2}\ead[label=e2]{LeonenkoN@cardiff.ac.uk}},
\author[C]{\fnms{Mar\'{i}a~D.}~\snm{Ruiz-Medina}\corref{}\thanksref{t2}\ead[label=e3]{mruiz@ugr.es}}
\and
\author[A]{\fnms{Irina N.} \snm{Savich}\ead[label=e4]{sim\_ka@i.ua}}
\thankstext{t2}{Supported in part by the European commission
PIRSES-GA-2008-230804
(Marie Curie), projects MTM2009-13393 of the DGI, and P09-FQM-5052 of
the Andalousian CICE, Spain, and the Australian
Research Council grants A10024117 and DP 0345577.}
\runauthor{Ivanov, Leonenko, Ruiz-Medina and Savich}
\affiliation{National Technical University of Ukraine,
Cardiff University, University~of~Granada~and National Technical
University~of~Ukraine}
\address[A]{A. V. Ivanov\\
I. N. Savich\\
National Technical University of Ukraine\\
Kyiv Polytechnic Institute\\
37 Peremogy Ave., Bldg 7, Room 435\\
03056 Kyiv\\
Ukraine\\
\printead{e1}\\
\hphantom{E-mail:\ }\printead*{e4}}

\address[B]{N. Leonenko\\
Cardiff School of Mathematics\\
Cardiff University\\
Senghennydd Road\\
Cardiff CF 24 4AG\\
United Kingdom\\
\printead{e2}}

\address[C]{M. D. Ruiz-Medina \\
Faculty of Sciences\\
University of Granada\\
Campus Fuente Nueva s/n\\
18071, Granada\\
Spain\\
\printead{e3}}

\end{aug}

\received{\smonth{11} \syear{2010}}
\revised{\smonth{4} \syear{2012}}

%
\begin{abstract}
The limit Gaussian distribution of multivariate weighted functionals
of nonlinear transformations of Gaussian stationary processes,
having multiple singular spectra, is derived, under very general
conditions on the weight function. This paper is motivated by its
potential applications in nonlinear regression, and asymptotic
inference on nonlinear functionals of Gaussian stationary processes
with singular spectra.\vspace*{-2pt}
\end{abstract}

%


\end{frontmatter}

\section{Introduction}\label{sec1}

During the last thirty years, a number of papers have been devoted
to limit theorems for nonlinear transformations of Gaussian
processes and random fields. The pioneering results are those of
Taqqu~\cite{r46,r47} and Dobrushin and Major~\cite{r15}, for
convergence to Gaussian and non-Gaussian distributions, under
long-range dependence, in terms of Hermite expansions, as well as Breuer
and Major~\cite{r11}, Ivanov and Leonenko~\cite{r27}, Chambers and
Slud~\cite{r13}, on convergence to the Gaussian distribution by
using diagram formulas or graphical methods. This line of research
continues to be of interest today; see Berman~\cite{r10} for
$m$-dependent approximation approach, Ho and Hsing~\cite{r23} for
martingale approach, Nualart and Peccati~\cite{r36} (see also
Peccati and Tudor~\cite{r40}) for the application of Malliavin
calculus, Nourdin and Peccati~\cite{r33} in relation to Stein's
method and exact Berry--Esseen asymptotics for functionals of
Gaussian fields, Avram, Leonenko and Sakhno~\cite{r8} for an
extension of graphical method for random fields, to name only a few
papers. The volume of Doukhan, Oppenheim and Taqqu~\cite{r17}
contains outstanding surveys of the field. Limit theorems for
weighted functionals of stochastic processes, and for processes with
seasonalities were considered by a number
of authors, including Rosenblatt~\cite{r42}, Oppenheim, Ould~Haye and Viano~\cite%
{r37}, Haye~\cite{r19}, and their references. Limit theorems for
nonlinear transformations of vector Gaussian processes have been
obtained by Arcones~\cite{r3}; see also his references.

In this paper, our main result (Theorem~\ref{thmain2b}) states the
convergence to the Gaussian distribution of the multivariate
weighted functionals of nonlinear transformations $\psi(\xi(t))$
of Gaussian stationary processes $\xi(t),$ with multiple
singularities in their spectra, having covariance function (c.f.)
belonging to a parametric family defined in Assumption (A2)
below. Here, $\psi\in L_{2}(\mathbb{R},\varphi(x)\,dx)$ [see Assumption~(A3)
in Section~\ref{sec2}], where $\varphi(x)=e^{-{x^{2}}/{2}%
}/\sqrt{2\pi},$ $x\in\mathbb{R},$ denotes the standard Gaussian
probability density. Specifically, under suitable conditions, the
convergence to the Gaussian distribution of
%
\begin{equation}
\bolds{\zeta}_{T}=W_{T}^{-1}\int
_{0}^{T}\mathbf{w}(t)\psi\bigl(\xi(t)\bigr)\nu
(dt)\label{W}
\end{equation}
as $T\rightarrow\infty,$ is obtained for certain ranges
of the
parameters defining the spectral singularities of $\xi$; see
Assumption (A4) in the next section. For each $T>0,$
%
\begin{eqnarray}\label{WF}
\mathbf{w}(t)&=& \bigl( w_{1}(t),\ldots,w_{q}(t)
\bigr)^{\prime},\qquad W_{T}^{2}=\operatorname{diag} \bigl(
W_{iT}^{2} \bigr)_{i=1}^{q},
\nonumber
\\[-8pt]
\\[-8pt]
\nonumber
W_{iT}^{2}&=&\int_{0}^{T}w_{i}^{2}(t)
\nu(dt),\qquad i=1,\ldots,q,
\end{eqnarray}
where, to ensure a finite limit variance, the weak
convergence of the family of matrix measures associated with
$\mathbf{w}$ over the intervals $\{[0,T], T>0\}$ is also assumed,
jointly with some restrictions on the boundedness of their
components in some neighborhoods of
the spectral singularities of $\xi$; see Assumptions (B1)
and (B2) in Section~\ref{sec4}. The convergence to the
Gaussian distribution also requires some conditions to be assumed on
the norms of
components of function $\mathbf{w}$; see condition (B3) in
Section %
\ref{sec5}.

As commented, the spectral density (s.d.) $f$ of $\xi(t)$ is
assumed to display several singularities denoted as $\Xi_{\mathrm{noise}}=\{\pm
\varkappa_{0},\ldots,\pm\varkappa_{r}\},$ with
$0\leq\varkappa_{0}<\varkappa_{1}< \cdots<\varkappa_{r}.$ In the
case where the weak-sense limit of the
measures associated with the multivariate weight function $\mathbf{w}
$ is an atomic measure, it is also assumed that its atoms $\Xi_{\mathrm{regr}}=\{
\delta_{1},\ldots,\delta_{n}\}$ do not intersect with
the singularities of $f$ (i.e., $\delta_{i}\neq\pm\varkappa_{j}, i=1,\ldots,n, j=0,1,\ldots,r$). The convergence to the
Gaussian distribution then holds with standard normalization.

The nature of the limit results obtained depends on the intersection
of the two spectral point sets $\Xi_{\mathrm{noise}}$ and $\Xi_{\mathrm{regr}}.$ In
the discrete case, this phenomenon was discussed by Yajima
\cite{r52,r53} in some other regression scheme. Otherwise, different
normalizing factors must be derived, and new limiting distributions
are obtained, for Hermite rank $m\geq2.$ Note that the classical
noncentral limit theorems (Taqqu~\cite{r47}, and Dobrushin and
Major~\cite{r15}) can be viewed as particular cases of the general
setting considered here, when there is an unique singular point in
the spectrum of $\xi,$ with $\varkappa_{0}=0,$ and $w(t)=1$. In
this case, the measure sequence, constructed from the weight
function $w,$ is given in terms of the Fejer kernel, which tends to
the delta-measure with atom at zero. Some limiting distributions for
the case when the two spectral point sets $\Xi_{\mathrm{noise}}$ and $\Xi
_{\mathrm{regr}}$ are in fact overlapped, in discrete time, can be derived
from the papers by Rosenblatt~\cite{r42}, Arcones~\cite{r3},
Oppenheim, Ould~Haye and Viano~\cite{r37} and Haye~\cite{r19}. In
continuous time, the limiting distributions for nonempty set, $\Xi
_{\mathrm{noise}}\cap\Xi_{\mathrm{regr}},$ can be obtained from the paper of
Ivanov and Leonenko~\cite{r29}, and the book by Ivanov and Leonenko
\cite{r27}. This subject will be considered in subsequent papers.

In the derivation of the main result of this paper, Peccati and Tudor's
central limit theorem~\cite{r40} (see also Nualart and Peccati \cite%
{r36}), for a family of vectors of random variables (r.v.'s)
belonging to fixed Wiener chaoses, is applied. The outline of the
paper is the following. Motivating examples, as well as preliminary
identities, and conditions needed in the derivation of the
subsequent results are provided in Section~\ref{sec2}. The zero-mean Gaussian
random field family considered is embedded into an isonormal process
family in Section~\ref{sec3}. The conditions needed for the weak-convergence
(in particular, to an atomic measure) of the matrix-valued measures
associated with the class of vectorial weight functions studied are
established in Section~\ref{sec4}. The asymptotic normality of the
corresponding weighted functionals of nonlinear transformations of
zero-mean Gaussian stationary random processes is obtained in
Section~\ref{sec5}. Section~\ref{sec6} provides the final comments, and our main
conjecture on the work is developed.

\section{Stationary processes with singular spectra}
\label{sec2} Let us consider simultaneously discrete and continuous
time
cases in the following development. Specifically, for a stationary
process $%
\xi$ defined on a complete probability space $(\Omega
,\mathbb{F},P),$ the following notation will be followed:
\[
\xi(t)=\xi(\omega,t)\dvtx\Omega\times\mathbb{S}\rightarrow\mathbb{R},
\]
where $\mathbb{S}=\mathbb{Z},$ for discrete time $t\in
\mathbb{Z},$ and $\mathbb{S}=\mathbb{R},$ for continuous time $t\in
\mathbb{R}.$ Such a process is assumed to be measurable and
mean-square continuous in the case of continuous time [see also
Assumption (A1) below].

In the definition of integrals, $\nu(dt)$ will represent a counting
measure in the case of discrete time [i.e., $\nu(\{t\})=1,$ $t\in
\mathbb{Z}$], and the Lebesgue measure $dt,$ in the~case of
continuous time [i.e., $\nu(dt)=dt $ if $t\in\mathbb{R}$].
According to this notation, the integral
\[
\int_{0}^{T}g(t)\xi(t)\nu(dt)
\]
represents the sum $\sum_{t=1}^{T}\xi(t)g(t),$ for
discrete time, and the Lebesgue integral $\int_{0}^{T}g(t)\xi
(t)\,dt,$ for continuous time, where $g(t)$ is a nonrandom
(measurable for continuous time) function.

Consider now the following motivating example.

\begin{example*}
Let $x$ be defined in terms of the
nonlinear regression model
%
\begin{equation}
x(t)=g(t,\theta)+\psi\bigl(\xi(t)\bigr), \qquad t \in\mathbb{S}_{+},
\label{Reg1}
\end{equation}
where $\mathbb{S}_{+} = \mathbb{R}_{+},$ for continuous
time and $\mathbb{S}_{+} = \mathbb{N},$ for
discrete time, and with $g(t,\theta)\dvtx\mathbb{S}_{+}\times\Theta
\rightarrow\mathbb{R}$
being a
continuously differentiable function of an unknown parameter $\theta
\in\Theta\subset\mathbb{R}^{q},$ consider $g_{i}(t,\theta
)=(\partial/\partial\theta_{i})g(t,\theta),$ $i=1,\ldots,q,$
such that
%
\begin{equation}
d_{iT}^{2}=\int_{0}^{T}
\bigl[ g_{i}(t,\theta) \bigr]^{2}\nu(dt)<\infty, \qquad T>0, i=1,
\ldots,q, \label{Reg2}
\end{equation}
and $\psi(\xi(t))$ represents the noise, with
$\mathit{E}\psi(\xi(t))=0.$ The least squares estimate (LSE)
$\hat{\theta}_{T}$ of an unknown parameter $\theta\in\Theta$,
obtained from the observations $x(t),$ $t \in[0,T],$ or
$t=1,\ldots,T,$ is any r.v. $\hat{\theta}_{T}\in\Theta^{c},$ having
the property
\[
Q_{T}(\hat{\theta}_{T})=\inf_{\tau\in\Theta^{c}}Q_{T}(
\tau),\qquad Q_{T}(\tau)=\int_{0}^{T}
\bigl[x(t)-g(t,\tau)\bigr]^{2}\nu(dt),
\]
where $\Theta^{c}$ is the closure of $\Theta$. Let $\nabla
g(t,\theta)= ( g_{1}(t,\theta),\ldots,g_{q}(t,\theta) )^{\prime}$ be the
column vector-gradient of the function $g(t,\theta).$ We denote $%
d_{T}^{2}(\theta)=\operatorname{diag} ( d_{iT}^{2} )_{i=1}^{q},$ where $%
d_{iT}^{2},$ $i=1,\ldots,q,$ are defined by (\ref{Reg2}). In the
theory of statistical estimation of unknown parameter $\theta\in
\Theta\subset
\mathbb{R}^{q}$ for the scheme (\ref{Reg1}), the asymptotic behavior,
as $%
T\rightarrow\infty,$ of the functional
%
\begin{equation}
\zeta_{T}=d_{T}^{-1}(\theta)\int
_{0}^{T}\nabla g(t,\theta)\psi\bigl(\xi(t)\bigr)
\nu(dt), \label{reg2}
\end{equation}
plays a crucial role, since, under certain number of
conditions,
the asymptotic distributions of the normalized LSE $d_{T}(\theta)(\hat{
\theta}_{T}-\theta),$ and properly normalized functional
(\ref{reg2}) coincide, as $T\rightarrow\infty$; see Ivanov and
Leonenko~\mbox{\cite{r27,r29}}.
\end{example*}

In this setting, an interesting case corresponds to $\xi(t)$ to be
a Gaussian stationary process with s.d. $f(\lambda)$ displaying
singularities
at the points $\Xi_{\mathrm{noise}}=\{\pm \varkappa_{j}, j=0,1,2,\ldots, r\}$; see
(%
\ref{sp1}) below. The nonlinear functions%
\[
g(t,\theta)=t^{\beta}\cos(t\vartheta+\phi),\qquad \beta\geq0, \vartheta
\in\mathbb{R}, \phi\in(-\pi,\pi], \theta=(\beta,\vartheta
,\phi)
\]
are of particular interest in applications because they
themselves also involve various seasonalities.

Let us consider $\{\xi(t), t\in\mathbb{S}\}$ to be a stochastic
process satisfying the following assumptions:
\begin{longlist}[(A1)]
\item[(A1)] Process $\xi$ is a real stationary mean-square
continuous Gaussian process with $\mathit{E}\xi(t)=0,$
$\mathit{E}\xi^{2}(t)=1.$

\item[(A2)] The c.f. of $\xi$ is of the form
%
\begin{equation}
B ( t ) =\mathit{E}\bigl[\xi(0)\xi(t)\bigr]=\sum
_{j=0}^{r}A_{j}B_{\alpha_{j},\varkappa_{j}}
( t ) ,\qquad  t\in\mathbb{S}, r\geq0, \label{cov}
\end{equation}
where, for $j=0,\ldots,r$,
\begin{eqnarray*}
B_{\alpha_{j},\varkappa_{j}} ( t ) &=&\frac{\cos (
\varkappa_{j}t ) }{ ( 1+t^{2} )^{\alpha
_{j}/2}}, \qquad 0\leq\varkappa_{0}<
\varkappa_{1}<\cdots<\varkappa_{r}, 0<\alpha_{j}<1,
t\in\mathbb{S},
\\
\sum_{j=0}^{r}A_{j}&=&1,\qquad
A_{j}\geq0, j=0,\ldots,r.
\end{eqnarray*}

The c.f. $B(t),$ $t\in\mathbb{S}$ admits the following spectral
decomposition:
\[
B(t)=\int_{\Lambda}e^{i\lambda t}f(\lambda)\,d\lambda,\qquad  t\in
\mathbb{R},
\]
where the set $\Lambda=(-\pi,\pi],$ in the discrete
case $(t\in\mathbb{Z}),$ and $\Lambda=\mathbb{R},$ in the
continuous case ($t\in\mathbb{R)}$, and the s.d. $f$ in the
continuous time is of the form
%
\begin{equation}
f ( \lambda) =\sum_{j=0}^{r}A_{j}f_{\alpha_{j},\varkappa
_{j}}
( \lambda) ,\qquad \lambda\in\mathbb{R}, \label{sp1}
\end{equation}
where, for $j=0,\ldots,r,$ and $\lambda\in\mathbb{R},$
\begin{eqnarray*}
f_{\alpha_{j},\varkappa_{j}} ( \lambda)& =&\frac{c_{1} ( \alpha_{j} )
}{2} \bigl[ K_{{(\alpha_{j}-1)}/{2}} \bigl(
\llvert\lambda+\varkappa_{j}\rrvert\bigr) \llvert\lambda+
\varkappa_{j}\rrvert^{{(\alpha
_{j}-1)}/{2}}\\
&&\hspace*{35pt}{}+K_{{(\alpha_{j}-1)}/{2}} \bigl( \llvert
\lambda-\varkappa_{j}\rrvert\bigr) \llvert\lambda-
\varkappa_{j}\rrvert^{{(\alpha_{j}-1)}/{2}} \bigr] ,
\end{eqnarray*}
with
\[
c_{1} ( \alpha) =\frac{2^{ ( 1-\alpha ) /2}}{\sqrt{\pi}%
\Gamma ( {\alpha}/{2} ) }
\]
and
\[
K_{\nu} ( z ) =\frac{1}{2}\int_{0}^{\infty}s^{\nu-1}
\exp\biggl\{ -\frac{1}{2} \biggl( s+\frac{1}{s} \biggr) z \biggr\}
\,ds,\qquad z\geq0, \nu\in\mathbb{R},
\]
being the modified Bessel function of the third kind of
order $\nu$ or McDonald's function. We also note that $K_{-\nu
} ( z ) =K_{\nu} ( z ) ,$ and for $z\downarrow0$
$K_{\nu} ( z ) \sim\Gamma ( \nu ) 2^{\nu
-1}z^{-\nu},$ $\nu>0.$

Thus, as $\lambda\rightarrow\pm\varkappa_{j},$ for $j=0,\ldots
,r,$
%
\begin{eqnarray}
f_{\alpha_{j},\varkappa_{j}} ( \lambda) &=& \frac{c_{2} ( \alpha_{j}
) }{2} \bigl[ \llvert\lambda+
\varkappa_{j}\rrvert^{\alpha_{j}-1} \bigl( 1-h_{j} \bigl(
\llvert\lambda+\varkappa_{j}\rrvert\bigr) \bigr)
\nonumber
\\[-8pt]
\\[-8pt]
\nonumber
&&\hspace*{36pt}{} +\llvert\lambda-\varkappa_{j}\rrvert^{\alpha_{j}-1}
\bigl( 1-h_{j} \bigl( \llvert\lambda-\varkappa_{j}\rrvert
\bigr) \bigr) \bigr],
\end{eqnarray}
where
\begin{eqnarray}
c_{2}(\alpha)&=&\frac{1}{2\Gamma(\alpha)\cos({\alpha\pi }/{2})},
\nonumber\\
h_{j} \bigl( \llvert\lambda\rrvert\bigr) &=&\frac{\Gamma
(
{(\alpha_{j}+1)}/{2} ) }{\Gamma ( {(3-\alpha_{j})}/{2}%
) }\biggl
\llvert\frac{\lambda}{2}\biggr\rrvert^{1-\alpha_{j}}+\frac{%
\Gamma ( {(\alpha_{j}+1)}/{2} ) }{4\Gamma (
{(3+\alpha_{j})}/{2} ) }\biggl
\llvert\frac{\lambda
}{2}\biggr\rrvert^{2}-o \bigl( \llvert\lambda
\rrvert^{2} \bigr) ,\nonumber\\
\eqntext{\lambda\longrightarrow0, j=0,\ldots,r.}
\end{eqnarray}
Therefore, the s.d. $f$ has $2r+2$ different singular points [see
condition (A2)], when $\varkappa_{0}\neq0.$

A model with discrete time which satisfies condition (A2) can
be obtained by using discretization procedure and the formula for
s.d. of
stationary processes with discrete time of the form%
\[
\sum_{k=-\infty}^{\infty}f ( \lambda+2k\pi) .
\]
 We will use the same notation for the s.d. in
both cases corresponding to discrete and continuous time.

Similar results can be obtained for c.f.'s of the form
\[
R_{\alpha_{j},\varkappa_{j}} ( t ) =\frac{\cos (
\varkappa_{j}t ) }{ ( 1+\llvert t\rrvert^{\rho
_{j}} )^{\alpha_{j}}}, \qquad\varkappa_{j}\in
\mathbb{R}, 0<\alpha_{j}\rho_{j}<1, \varkappa_{j}
\neq0, j=0,\ldots,r
\]
(see again Ivanov and Leonenko~\cite{r29} for details).

It is well known that the Hermite polynomials $H_{k}(x)=(-1)^{k}e^{
{%
x^{2}}/{2}}\frac{d^{k}}{dx^{k}}e^{-{x^{2}}/{2}},$ $k=0,1,\ldots$
constitute a complete orthogonal system in the Hilbert space
$L_{2}(\mathbb{R},\varphi(x)\,dx)$ of square integrable functions with respect to the
standard Gaussian density $\varphi.$

\item[(A3)] Assume that the function $\psi\in
L_{2}(\mathbb{R},\varphi
(x)\,dx)$, that is, $\mathit{E}\psi^{2}(\xi(0))<\infty,$ and
$C_{0}(\psi)=%
\mathit{E}\psi(\xi(0))=0.$

%
\begin{definition}
A function $\psi\in L_{2}(\mathbb{R},\varphi(x)\,dx)$ has Hermite rank\break
$H\operatorname{rank}(\psi)=m$ if either $C_{1}(\psi)\neq0$ and $m=1,$ or for some $%
m\geq2,$ $C_{1}(\psi)=\cdots=C_{m-1}(\psi)=0, C_{m}(\psi)\neq
0.$
\end{definition}

\item[(A4)] Either (i) $H\operatorname{rank}(\psi)=1,$ $\alpha>1/2;$ or (ii)
$H\operatorname{rank}(\psi)=m,$ $\alpha m>1,$ where $\alpha
=\min_{j=0,\ldots
,r}\alpha_{j},$ with $\alpha_{j},$ $j=0,\ldots,r,$ introduced in
(A2).

Under condition (A3), function $\psi(x)$ of $H\operatorname{rank}(\psi
)=m$ can
be expanded into a Hermite series in the Hilbert space $L_{2}(\mathbb{R}
,\varphi(x)\,dx)$
%
\begin{equation}
\psi(x)=\sum_{k=m}^{\infty}
\frac{C_{k}(\psi)}{k!}H_{k}(x), \label{ser1}
\end{equation}
or the process $\psi(\xi(t))$ admits a Hermite series expansion in
the Hilbert space $L_{2}(\Omega,\mathbb{F},P)$
%
\begin{equation}
\psi\bigl(\xi(t)\bigr)=\sum_{k=m}^{\infty}
\frac{C_{k}(\psi
)}{k!}H_{k}\bigl(\xi(t)\bigr), \label{ser2}
\end{equation}
where
\[
C_{k}(\psi)=\int_{\mathbb{R}}\psi(x)H_{k}(x)
\varphi(x)\,dx,\qquad  k\geq0.
\]
\end{longlist}
%
\section{Some elements of the theory of isonormal processes}\label{sec3}

In this section, we introduce basic notation, elements and results
in relation to Gaussian Hilbert spaces, isonormal processes and
chaos expansions needed for our purposes; see Nualart~\cite{r35};
Janson~\cite{r30}; Nualart and Peccati~\cite{r36}; Peccati and Tudor
\cite{r40}; Peccati~\cite{r38}; Nourdin, Peccati and R{\'e}veillac
\cite{r34}, among others.

%
\begin{definition}
Let $H$ be a real separable Hilbert space. The set of r.v.'s
$X=\{X(h)\dvtx h\in
H\}$ is said to be an isonormal process on $H$ if $X$ is a centered $H$-indexed
Gaussian family defined on a probability space $(\Omega,\mathbb
{F},%
\mathbb{P}),$ and it satisfies
\[
\mathit{E}\bigl[X(h)X(g)\bigr]= \langle h,g \rangle_{H},\qquad h,g\in H.
\]
\end{definition}

Let us now consider a real-valued centered Gaussian process $\xi$
indexed over $\mathbb{S}=\mathbb{R}.$ By $\mathcal{E}$ denote the
collection of all
finite linear combinations of indicator functions of the type $\mathbf
{l}%
_{(-\infty,t]},$ with $t\in\mathbb{R}.$ To embed a real-valued
centered
Gaussian process $\xi$ indexed by $\mathbb{R}$ into some isonormal
process $%
X,$ we introduce a separable Hilbert space $H$ defined as the closure
of $%
\mathcal{E}$ with respect to the scalar product
%
\begin{equation}
\langle f,h \rangle_{H}:=\sum_{i,j}a_{i}c_{j}
\mathit{E}\xi(s_{i})\xi(t_{j}) \label{innisonormalproc}
\end{equation}
 for given functions $f=\sum_{i}a_{i}\mathbf{l}_{(-\infty
,s_{i}]} $ and $h=\sum_{j}c_{j}\mathbf{l}_{(-\infty,t_{j}]}$ in
$H.$ Thus,
for any function $h=\sum_{i}c_{i}\mathbf{l}_{(-\infty,t_{i}]}\in
\mathcal{E}%
,$ define
%
\begin{equation}
X(h)=\sum_{i}c_{i}\xi
(t_{i}). \label{isoproc}
\end{equation}

Additionally, for any function $h\in H,$ $X(h)$ can be defined as
the limit
in $L_{2}(\Omega,\mathbb{F},\mathbb{P})$ of $X(h_{n})$ for any
sequence\vadjust{\goodbreak} $%
\{h_{n}\}\subset\mathcal{E}$ convergent to $h$ in $H.$ This
sequence may be not unique, but the definition of $X(h)$ does not
depend on the choice of the sequence $\{h_{n}\}.$ From this
construction, process $X$ is an isonormal process over $H$ defined
as $X(\mathbf{l}_{(-\infty,t]})=\xi(t).$

When $\mathbb{S}=\mathbb{Z},$ a similar development in terms of
sequences leads to the definition of an isonormal process from a
Gaussian process $\xi$ on $\mathbb{Z}.$ Now, $\mathcal{E}$ denotes
the set of all real-valued sequences $h=\{h_{l}\dvtx l\in\mathbb{Z}\}$
such that $h_{l}\neq0$ only for a finite number of integers $l.$
The real separable Hilbert space $H$ is then introduced as the
closure of the set $\mathcal{E}$ with respect to the scalar product
\[
\langle f,h \rangle_{H}:=\sum_{k,l}f_{k}h_{l}
\mathit{E}\xi(k)\xi(l)
\]
 for given sequences $f=\{f_{k}\dvtx k\in\mathbb{Z}\}$ and
$h=\{h_{l}\dvtx l\in\mathbb{Z}\}.$
If $h\in H,$ then the series $\sum_{l\in\mathbb{Z}}h_{l}\xi(l)$
converges in $L_{2}(\Omega,\mathbb{F},\mathbb{P}).$ Thus the
centered Gaussian family $\{X(h)\dvtx h\in H\},$ with
%
\begin{equation}
X(h)=\sum_{l\in\mathbb{Z}}h_{l}\xi(l) \label{iso2}
\end{equation}
 is an isonormal process over $H.$

Let $X$ be an isonormal process defined on $H$ as before, that is, from
a
centered Gaussian random process $\xi.$ Let us write $\mathcal
{H}_{0}(X)=%
\mathbb{R}^{1},$ and $\mathcal{H}_{1}(X)$ the closed linear
subspace of the
set of r.v.'s $\{X(h)\dvtx h\in H\}$ in the Hilbert space $L_{2}(\Omega
,\mathbb{F%
},\mathbb{P}).$ Thus
\begin{eqnarray*}
X\dvtx H&\longrightarrow&\mathcal{H}_{1}(X),
\\
h&\longrightarrow& X(h).
\end{eqnarray*}

For any $n\geq2,$ by $\mathcal{H}_{n}(X),$ the $n$th Wiener chaos
of process $X$ is denoted, that is, the closed subspace of
$L_{2}(\Omega,%
\mathbb{F},\mathbb{P})$ generated by the r.v.'s $H_{n}(Y),$ where
$Y\in\mathcal{H}_{1}(X),$ and $\mathit{E}[Y^{2}]=1,$ with $H_{n}$
denoting, as before, the $n$th Hermite polynomial. Let us now
consider the isometry
%
\begin{equation}
I_{n}^{X}\dvtx H^{\odot n}\longrightarrow
\mathcal{H}_{n}(X), \label{isometry}
\end{equation}
between the symmetric tensor product $H^{\odot n},$
equipped with the norm $\sqrt{n!}\Vert\cdot\Vert
_{H^{\otimes n}},$ and the $n$th Wiener
chaos $\mathcal{H}_{n}(X)$ of $X.$ For any $h\in H^{\otimes n},$ $%
I_{n}^{X}(h)$ is then defined as
$I_{n}^{X}(h):=I_{n}^{X}(\widetilde{h}),$ with $\widetilde{h}$
denoting the symmetrization of $h.$ For any $g\in H^{\otimes m}$ and
$h\in H^{\otimes n},$
\[
E\bigl[I_{m}^{X}(g)I_{n}^{X}(h)
\bigr]=\delta_{mn}m! \langle\widetilde{g},%
\widetilde{h}
\rangle_{H^{\otimes m}}.
\]

The $p$th contraction of $g=g_{1}\otimes\cdots\otimes g_{k}\in
H^{\otimes k}$ and $h=h_{1}\otimes\cdots\otimes h_{k}\in H^{\otimes
k},$ designated as $g\otimes_{p}h,$ is the element of $H^{\otimes
2(k-p)}$ given by
%
\begin{equation}
\qquad g\otimes_{p}h= \langle h_{1},g_{1}
\rangle_{H}\cdots\langle h_{p},g_{p}
\rangle_{H}g_{p+1}\otimes\cdots\otimes g_{k}
\otimes h_{p+1}\otimes\cdots\otimes h_{k}. \label{eqcontp}
\end{equation}

The definition can be extended by linearity to any element of $%
H^{\otimes k}.$ Finally, any r.v. $%
F\in L_{2}(\Omega,\mathbb{G},\mathbb{P}),$ with $\sigma$- field
$\mathbb{G}$ generated by the r.v.'s $\{X(h),h\in H\},$ admits an
unique chaos decomposition $F=\sum_{k=0}^{\infty}I_{k}^{X}(h_{k}),$
where \mbox{$h_{k}\in H^{\odot k}$}.\vadjust{\goodbreak}

From the constructions (\ref{isoproc}) and (\ref{iso2}) of an
isonormal process $X$ from a Gaussian process $\xi$, respectively,
defined over continuous and discrete time,
$\mathcal{H}_{n}(X), n\geq1,$ coincides with the $n$th Wiener chaos
associated with $\xi,$ $\mathcal{H}_{n}(\xi), n\geq1.$
Since, by definitions (\ref{isoproc}) and (\ref{iso2}), $\mathcal
{H}_{1}(X)=\mathcal{H}%
_{1}(\xi),$ and, as stated before, the $n$th Wiener chaos of
process $X$ is the closed subspace of $L_{2}(\Omega
,\mathbb{F},\mathbb{P})$ generated from the evaluation of $n$th
Hermite polynomial $H_{n}$ over the r.v.'s of the space
$\mathcal{H}_{1}(X)=\mathcal{H}_{1}(\xi).$

The next statement is a convenient, for our purposes, modification
of Theorem~1 of Peccati and Tudor~\cite{r40}; see also Nualart and
Peccati~\cite{r36} (in the above papers all statements are
formulated for positive integers $T\in\{1,2,\ldots\}$, but it is easy
to see that one can formulate similar results for continuous $T>0$
as well).

%
\begin{proposition}
\label{pr2} Let $ \{ \xi(t), t\in\mathbb{S} \} $ be a
centered Gaussian process, and $X$ is the isonormal process
constructed from it as given in (\ref{isoproc}) and (\ref{iso2}).
Consider the natural numbers:
$1\leq n_{1}<n_{2}<\cdots<n_{d}<\infty, d\geq2,$ and the set of
r.v.'s $%
\pi_{T,n_{j}}(\xi)\in\mathcal{H}_{n_{j}}(\xi),$ where, for
$T>0,$ $\pi_{T,n_{j}}(\xi)=I_{n_{j}}^{X}(f_{j,T}),$ for certain
$f_{j,T}\in H^{\odot n_{j}},$ $j=1,\ldots,d,$ such that
%
\begin{equation}
\lim_{T\rightarrow\infty}\mathit{E}\pi_{T,n_{j}}^{2}(\xi)=
\lim_{T\rightarrow\infty}n_{j}!\Vert f_{j,T}\Vert_{H^{\otimes
n_{j}}}^{2}=1,\qquad
j=1,\ldots,d. \label{NP1}
\end{equation}

Then the following conditions are equivalent:

\begin{longlist}[(iii)]
\item[(i)] For each $j=1,\ldots,d,$
\[
\lim_{T\rightarrow\infty}\| f_{j,T}\otimes_{p} f_{j,T}
\|_{H^{\otimes2(n_{j}-p)}}=0
\]
 for every $p=1,\ldots,n_{j}-1.$

\item[(ii)] For every $j=1,\ldots,d,$
\[
\lim_{T\rightarrow\infty}E \bigl[ \bigl(I_{n_{j}}^{X}(f_{j,T})
\bigr)^{4}%
\bigr]=3.
\]

\item[(iii)] As $T\rightarrow\infty,$ the vector $%
(I_{n_{1}}^{X}(f_{1,T}),\ldots, I^{X}_{n_{d}}(f_{d,T}) )$
converges
in distribution to a $d$-dimensional standard Gaussian vector $N_{d}(0,%
\mathbb{I}_{d}).$
\end{longlist}
\end{proposition}

The proof follows from Peccati and Tudor~\cite{r40}, and
Nualart and Peccati~\cite{r36}, considering the fact that
$\mathcal{H}_{n}(\xi)=\mathcal{H}_{n}(X),$ for any $n\geq1,$
with $X$ being the isonormal process constructed from identity (\ref
{isoproc}%
), in the continuous time case, and, similarly, in the discrete time
case, from equation (\ref{iso2}).

%
\begin{corollary}
\label{cor1} Assume that conditions (\ref{NP1}) and \textup{(i)} or \textup{(ii)} of
Proposition~\ref{pr2} are satisfied for r.v.'s
%
\begin{equation}
\pi_{T,n_{j}}(\xi)=\int_{0}^{T}r_{T,j}(t)H_{n_{j}}
\bigl(\xi(t)\bigr)\nu(dt), \label{NP8}
\end{equation}
where, in the case of continuous time, it is also assumed
that $r_{T,j}(t)\in C([0,\infty)),$ for $T>0,$ and $j=1,\ldots,d.$
Then, the vector
%
\begin{equation}
\qquad \pi_{T,d}(\xi)= \biggl( \int_{0}^{T}r_{T,1}(t)H_{n_{1}}
\bigl(\xi(t)\bigr)\nu(dt),\ldots,\int_{0}^{T}r_{T,d}(t)H_{n_{d}}
\bigl(\xi(t)\bigr)\nu(dt) \biggr) \label{NP9}
\end{equation}
converges in distributions, as $T\rightarrow\infty,$ to a standard
Gaussian vector $\pi_{d}\sim N(0,\mathbb{I}_{d}).$
\end{corollary}

\begin{pf}
In the case of continuous time, since $\xi(t)=X ( \mathbf{l}%
_{(-\infty,t]} ) ,$
\[
H_{n_{j}}\bigl(\xi
(t)\bigr)=H_{n_{j}} \bigl(X ( \mathbf{l}%
_{(-\infty,t]} ) \bigr)= I_{n_{j}}^{X} \bigl(
\mathbf{l}_{(-\infty,t]}^{\otimes n_{j}} \bigr),
\]
 where
$I_{n_{j}}^{X}$ denotes the isometry introduced in (\ref{isometry}).
Therefore, for $r_{T,j}(t)\in C([0,\infty)),$ $T>0$ and for
$j=1,\ldots,d,$
\[
\pi_{T,n_{j}}(\xi)=\int_{0}^{T}r_{T,j}(t)H_{n_{j}}
\bigl(\xi(t)\bigr)\,dt=I_{n_{j}}^{X} \biggl( \int
_{0}^{T}r_{T,j}(t)\mathbf{l}_{(-\infty,t]}^{\otimes
n_{j}}\,dt
\biggr) . \label{ID}
\]
Thus, considering in (iii) of Proposition~\ref{pr2}
\[
f_{j,T}(s_{1},\ldots,s_{n_{j}})=\int_{0}^{T}r_{T,j}(t)\mathbf
{l}_{(-\infty,t]}^{\otimes
n_{j}}(s_{1},\ldots,s_{n_{j}})\,dt
\]
 for $j=1,\ldots,d,$ we
obtain the desired result.

Similarly, for the case of discrete time, we have, from
(\ref{iso2}),
%
\begin{equation}
X(\delta_{\bolds{\cdot},l})=\xi(l),\qquad  l\in\mathbb{Z}, \label{dti}
\end{equation}
 where, for each $l\in\mathbb{Z},$
$\delta_{\bolds{\cdot},l}$ denotes the Kronecker delta
function, that is,
\[
\delta_{i,l}= \cases{
1, &\quad$\mbox{if } i=l,$\vspace*{2pt}\cr
0, &\quad $\mbox{if } i\neq l, i \in\mathbb{Z}.$}
\]

Therefore,
\[
H_{n_{j}}\bigl(\xi(l)\bigr)=H_{n_{j}}\bigl (X (
\delta_{\bolds{\cdot},l} ) \bigr)= I_{n_{j}}^{X} \bigl(
\delta_{\bolds{\cdot},l}^{\otimes n_{j}} \bigr).
\]

Proposition~\ref{pr2}(iii) is then applied, considering
\[
f_{j,T}(m_{1},\ldots,m_{n_{j}})=\sum_{l=1}^{T}r_{T,j}(l)\prod
_{i=1}^{j}\delta_{m_{i},l}, \qquad m_{1},\ldots,m_{j}\in
\mathbb{Z}
\]
 for $j=1,\ldots,d.$
\end{pf}

\section{Spectral measures of weight functions and admissible spectral
densities}
\label{sec4}

Let us first establish some results on weak-convergence
of
matrix-valued measures, given by
%
\begin{equation}
\mu_{T}^{jl}(d\lambda)=\frac{w_{T}^{j}(\lambda)\overline{%
w_{T}^{l}(\lambda)}\,d\lambda}{\sqrt{\int_{\Lambda
}\llvert w_{T}^{j}(\lambda)\rrvert^{2}\,d\lambda
\int_{\Lambda}\llvert w_{T}^{l}(\lambda)\rrvert^{2}\,d\lambda}},\qquad
j,l=1,\ldots,q,
\label{32}\vadjust{\goodbreak}
\end{equation}
 where
\[
w_{T}^{j}(\lambda)=\int_{0}^{T}e^{\mathrm{i}t\lambda
}w_{j}(t)
\nu(dt),\qquad j=1,\ldots,q,
\]
and the functions $w_{j}(t),$ $j=1,\ldots,q,$ are, as before, the
functions (%
\ref{WF}) involved in the definition of the random vector (\ref{W}).

(B1) Assume that the weak-convergence $\bolds{\mu}%
_{T}\Rightarrow\bolds{\mu},$ when $T\rightarrow\infty$
holds,
where $\bolds{\mu}_{T}$ is defined by (\ref{32}) and $\bolds{
\mu}$ is a positive definite matrix measure.

The above condition means that an element $\mu^{jl}$ of the
matrix-valued measure $\bolds{\mu}$ is a signed measure of
bounded variation, and the matrix $\bolds{\mu}(A)$ is positive
definite for any set $A\in\mathcal{A},$ with $\mathcal{A}$ denoting
the $\sigma$-algebra of measurable subsets of $\mathbb{R}$; see, for
example, Ibragimov and Rozanov~\cite{r24}.

The following definition can be found in Grenander and Rosenblatt \cite
{r18}%
, Ibragimov and Rozanov~\cite{r24} and Ivanov and Leonenko
\cite{r27}.

%
\begin{definition}
\label{mvmeasure} The nondegenerate matrix-valued measure $\bolds
{\mu}%
(d\lambda)= \{\mu^{jl}(d\lambda)\}_{j,l=1}^{q}$ is said to be the
spectral measure of function $\mathbf{w(}t).$
\end{definition}

%
\begin{definition}[(Ibragimov and Rozanov~\cite{r24})] The s.d. $f$ is said to be
$\bolds{\mu}$-admissible if it is integrable, that is, all elements of the
matrix%
\[
\int_{\Lambda}f(\lambda)\bolds{\mu}(d\lambda)
\]
 are finite, and
%
\begin{equation}
\lim_{T\rightarrow\infty}\int_{\Lambda}f(\lambda)\bolds{\mu
}%
_{T}(d\lambda)=\int_{\Lambda}f(\lambda)
\bolds{\mu}(d\lambda). \label{star}
\end{equation}
\end{definition}

Let us introduce two conditions on the s.d. $f$ that guarantee its $%
\bolds{\mu}$-admissibility. These assumptions are related to
basic conditions on the c.f. and s.d. (A2). In the following,
$J$ denotes
one of the three sets:%
\[
\{-r,\ldots,-1,0,1,\ldots,r\};\{-r,\ldots,-1,1,\ldots,r\};\{0\}.
\]
We formulate the following condition for a set $J=\{-r,\ldots
,-1,0,1,\ldots,r%
\}.$

(I) The s.d. $f\in\mathcal{C} ( \Lambda\setminus
\{\varkappa_{j},j\in J\} ) ,$ with
\[
\varkappa_{-j}=-\varkappa_{j},\qquad j=0,1,\ldots,r, 0\leq
\varkappa_{0}<\varkappa_{1}<\cdots<\varkappa_{r}
\]
 and, for $j=0,1,\ldots,r,$
%
\begin{eqnarray}\label{sp3}
\lim_{\lambda\rightarrow\varkappa_{j}}f(\lambda)\llvert\lambda
-\varkappa_{j}
\rrvert^{1-\alpha_{j}}=a_{j}>0,
\nonumber
\\[-8pt]
\\[-8pt]
\eqntext{\alpha_{j}\in(0,1), j\in
J; \alpha_{-j}=\alpha_{j}, a_{-j}=a_{j}.}
\end{eqnarray}

We obtain from (\ref{sp3}) that, for any $\varepsilon>0,$ and $j\in
J,$
there exists $\delta_{j}=\delta_{j}(\varepsilon),$ such that for $%
\llvert \lambda-\varkappa_{j}\rrvert <\delta_{j}$
\[
f(\lambda)<\frac{a_{j}+\varepsilon}{\llvert \lambda-\varkappa
_{j}\rrvert^{1-\alpha_{j}}}.\vadjust{\goodbreak}
\]

Then, for $\llvert \lambda-\varkappa_{j}\rrvert <\delta_{j},$ we
have the following:
\[
\bigl\{\lambda\dvtx f(\lambda)>c\bigr\}\subset V_{j}(c)= \biggl\{
\lambda
\dvtx\llvert\lambda-\varkappa_{j}\rrvert< \biggl( \frac
{a_{j}+\varepsilon}{c}%
\biggr)^{{1}/{(1-\alpha_{j})}} \biggr\}.
\]

Moreover $c$ must satisfy the inequality%
\[
\biggl( \frac{a_{j}+\varepsilon}{c} \biggr)^{{1}/{(1-\alpha
_{j})}}\leq\delta_{j},
\]
and equivalently,
%
\begin{equation}
c\geq\frac{a_{j}+\varepsilon}{\delta_{j}^{1-\alpha_{j}}(\varepsilon
)}%
=c_{j}(\varepsilon). \label{sp4}
\end{equation}

(II) Let $\varepsilon_{0}>0$ be fixed. There exists $%
c_{0}=\max_{j\in J}c_{j}(\varepsilon_{0}),$ such that for $c\geq c_{0},$
%
\begin{equation}
\bigl\{\lambda\dvtx f(\lambda)>c\bigr\}\subset\bigcup_{j\in J}V_{j}(c),
\label{sp5}
\end{equation}
where $c_{j}(\varepsilon)$ are defined by (\ref{sp4}).

It is easy to see that for sufficiently large $c$ (say, $c\geq
c_{0}),$ the neighborhoods $V_{j}(c),$ $j\in J,$ in (\ref{sp5})$,$
are nonoverlapping, and
\[
\bigl\llvert V_{j}(c)\bigr\rrvert\downarrow0
\]
as $c\rightarrow\infty.$

Note that the function (\ref{sp1}) satisfies conditions
(I) and (II).

(B2) For $T$ sufficiently large (say, $T\geq
T_{0}),$
%
\begin{equation}
W_{iT}^{-1}\max_{\lambda\in V_{j}(c_{0})}\bigl\llvert
w_{T}^{i}(\lambda)\bigr\rrvert\leq k_{ij}<
\infty,\qquad j\in J, i=1,\ldots,q. \label{sp6}
\end{equation}

In condition (B2), one can assume that (\ref{sp6}) holds only
for $%
j=0,1,\ldots,r,$ since $V_{-j}(c_{0})=-V_{j}(c_{0}),$ $j=0,1,\ldots
,r.$

%
\begin{theorem}
\label{atresul} Assume that conditions \textup{(B1)}, \textup{(B2)}, as well as
\textup{(I)}, \textup{(II)} are satisfied, and the s.d. $f$ is integrable with respect to
the
spectral measure $\bolds{\mu}$, then the s.d. $f$ is
$\bolds{\mu}$-admissible.
\end{theorem}

\begin{pf} For $c\geq c_{0},$ we consider
\[
f^{c}(\lambda)=f(\lambda){1}_{ \langle f(\lambda)< c
\rangle}(\lambda)+c{
1}_{ \langle f(\lambda)\geq
c \rangle}(\lambda).
\]

Then, for $k,l=1,\ldots,q,$
%
\begin{eqnarray}\label{eqintidb}
&&\biggl\llvert\int_{\Lambda}f(\lambda)\mu_{T}^{k,l}(d
\lambda)-\int_{\Lambda}f(\lambda)\mu^{k,l}(d\lambda)
\biggr\rrvert
\nonumber
\\
&&\qquad\leq\biggl\llvert\int_{\Lambda}f(\lambda)
\mu_{T}^{k,l}(d\lambda)-\int_{\Lambda}f^{c}(
\lambda)\mu_{T}^{k,l}(d\lambda)\biggr\rrvert
\nonumber
\\
&&\qquad\quad{}+\biggl\llvert\int_{\Lambda}f^{c}(\lambda)
\mu_{T}^{k,l}(d\lambda)-\int_{\Lambda}f^{c}(
\lambda)\mu^{k,l}(d\lambda)\biggr\rrvert
\\
&&\qquad\quad{}+\biggl\llvert\int_{\Lambda}f^{c}(\lambda)
\mu^{k,l}(d\lambda)-\int_{\Lambda}f(\lambda)
\mu^{k,l}(d\lambda)\biggr\rrvert\nonumber\\
&&\qquad=I_{1}^{k,l}(T,c)+I_{2}^{k,l}(T,c)+I_{3}^{k,l}(c).\nonumber
\end{eqnarray}
By Assumption (B1), for any complex numbers
$z=(z_{1},\ldots,z_{q}),$
the function%
\[
M_{z}(A)=\sum_{k,l=1}^{q}
\mu^{k,l}(A)z_{k}\bar{z}_{l}\geq0,\qquad A\in
\mathcal{A},
\]
 is a measure. Thus, by Lebesgue's monotone convergence
theorem,
%
\begin{equation}
\int_{\Lambda}f^{c}(\lambda)M_{z}(d\lambda
)\mathop{\longrightarrow}_{c\rightarrow
\infty}\int_{\Lambda}f(\lambda
)M_{z}(d\lambda). \label{sp8}
\end{equation}
Note that the diagonal elements $\mu^{k,k}$ and $\mu^{l,l}$ are
measures; thus if only $z_{k}$ and $z_{l}$ are nonzero among
$z=(z_{1},\ldots
,z_{q}),$ we obtain from (\ref{sp8}) that%
%
\begin{equation}\qquad
\int_{\Lambda} \bigl( f(\lambda)-f^{c}(\lambda)
\bigr) \bigl( \mu^{k,l}(d\lambda)z_{k}\bar{z}_{l}+
\mu^{l,k}(d\lambda)z_{l}\bar{z}%
_{k}
\bigr) \longrightarrow0,\qquad c\rightarrow\infty. \label{sp9}
\end{equation}
Note that $\mu^{l,k}=\overline{\mu^{k,l}},$ and choosing, for instance,
$%
z_{k}=z_{l}=1,$ we have from (\ref{sp9}),%
\[
\int_{\Lambda} \bigl( f(\lambda)-f^{c}(\lambda)
\bigr) \operatorname{Re}\bigl(\mu^{k,l}\bigr) (d\lambda)\longrightarrow
0,\qquad c\rightarrow\infty.
\]
If we choose $z_{k}=1, z_{l}=-\mathrm{i},$ then
\[
\int_{\Lambda} \bigl( f(\lambda)-f^{c}(\lambda)
\bigr) \operatorname{Im}\bigl(\mu^{k,l}\bigr) (d\lambda)\longrightarrow
0, \qquad c\rightarrow\infty.
\]
Thus%
\[
\lim_{c\rightarrow\infty}I_{3}^{k,l}(c)=0.
\]
For a fixed $c,$ we obtain, from condition (B1), that
\[
\lim_{T\rightarrow\infty}I_{2}^{k,l}(T,c)=0.
\]

On the other hand, under the conditions assumed in this theorem, for
$T\geq T_{0},$
%
\begin{eqnarray}
 I_{1}^{k,l}(T,c)&\leq&\frac{1}{2\pi}\int
_{ \{ \lambda
\dvtx f(\lambda)>c \} }\bigl(f(\lambda)-c\bigr)\frac{\llvert
w_{T}^{k}(\lambda)\rrvert \llvert w_{T}^{l}(\lambda
)\rrvert }{W_{k,T} W_{l,T}}\,d\lambda
\nonumber
\\[-8pt]
\\[-8pt]
\nonumber
& \leq&\frac{1}{2\pi}\sum_{j\in
J}k_{j,k}k_{j,l}
\int_{V_{j}(c)}f(\lambda)\,d\lambda\rightarrow0,
\end{eqnarray}
when $c\rightarrow\infty.$ Thus, for any $\varepsilon>0$
and $T\geq T_{0},$ one can choose $%
c_{1}=c_{1}(\varepsilon)\geq c_{0},$ such that for $c>c_{1},$ we have $
I_{1}^{k,l}(T,c)<\varepsilon/3.$ Then, once can take
$c_{2}=c_{2}(\varepsilon)\geq c_{0},$ such that for $%
c>c_{2},$ we have $I_{3}^{k,l}(c)<\varepsilon/3.$

Let us now fix $c=\max(c_{1},c_{2});$ then, there exists $T_{1}=T_{1}(%
\varepsilon)>T_{0},$ such that for $T>T_{0},$
$I_{2}^{k,l}(T,c)<\varepsilon/3,$ and the left-hand side of
(\ref{eqintidb}) is less than~$\varepsilon.$
\end{pf}

\section{Central limit theorem for weighted functionals}
\label{sec5}

This section provides the asymptotic normality as
$T\rightarrow
\infty$ of the vector (\ref{W}), that is, we will prove that the
vector $%
\bolds{\zeta}_{T}$ converges in distribution ($\Longrightarrow$)
to some Gaussian vector $\bolds{\zeta}.$ Thus, for any $z\in
\mathbb{R}^{q},$ we prove that $ \langle\bolds{\zeta
}_{T},z \rangle
\Longrightarrow \langle\bolds{\zeta},z \rangle,$ as $%
T\rightarrow\infty.$ Denoting, for $z=(z_{1},\ldots,z_{q}),$
\[
\sum_{i=1}^{q}z_{i}W_{iT}^{-1}w_{i}(t)=R_{T}(t,z)=R_{T}(t),
\]
from (\ref{ser2}), we have
\[
\langle\bolds{\zeta}_{T},z \rangle=\int_{0}^{T}
\psi\bigl(\xi(t)\bigr)R_{T}(t)\nu(dt)=\sum
_{j=m}^{\infty}\frac{C_{j}(\psi)}{j!}%
\int
_{0}^{T}R_{T}(t)H_{j}\bigl(\xi
(t)\bigr)\nu(dt).
\]

In the derivation of the proof of our main result, the following
additional conditions are required:

(B3) For $T>T_{0},$%
%
\begin{equation}
W_{i,T}^{-1}\sup\bigl\llvert w_{i}(t)\bigr\rrvert
\leq k_{i}T^{-1/2},\qquad i=1,\ldots,q, \label{lim8}
\end{equation}
 where the supremum is taken over $t$ in the interval
$[0,T],$ in the case of continuous time, and over $t$ in the set
$\{1,\ldots,T\},$ in the case of discrete time.

Let $f^{(*1)}(\lambda)=f(\lambda),$ and for $j\geq2,$
\[
f^{*(j)}(\lambda)=\int_{\Lambda^{j-1}}f(\lambda-
\lambda_{2}-\cdots-\lambda_{j})\prod
_{i=2}^{j}f(\lambda_{i})\,d
\lambda_{2}%
\cdots d\lambda_{j},
\]
 the $j$th convolution of the s.d. $f(\lambda).$

(C) The matrix integrals
\[
\int_{\Lambda}f^{\ast(j)}(\lambda)\bolds{\mu}(d
\lambda),\qquad j\geq1,
\]
 are positive definite.

We now proceed the formulation of our main result.

%
\begin{theorem}
\label{thmain2b} Suppose that conditions \textup{(A1)--(A4)},
\textup{(B1)--(B3)} and \textup{(C)} are fulfilled. Then, the r.v.
$\bolds{\zeta}_{T}$ in (\ref{W}) converges in distribution, as
$T\longrightarrow\infty,$ to the Gaussian r.v. $\bolds{\zeta}$
with zero mean and covariance matrix
%
\begin{equation}
\Xi=2\pi\sum_{j=m}^{\infty}
\frac{C_{j}^{2}(\psi)}{j!}%
\int_{\Lambda}f^{\ast(j)}(
\lambda)\bolds{\mu}(d\lambda), \label{lim14}
\end{equation}
 where $\bolds{\mu}$ is the weak-sense limit of the
family of matrix-valued measures introduced in
(\ref{32}) and associated with the weight function $\mathbf{w}(t)$ in
(\ref{WF}%
), given from functional~(\ref{W}).
\end{theorem}

In the proof of the above theorem, the following identities will be
applied jointly with Lemma~\ref{lemprevth} formulated below.
Specifically, from the orthogonality of Hermite polynomials, we
obtain
%
\begin{eqnarray}\label{lim3}
\mathit{E} \langle\bolds{\zeta}_{T},z \rangle^{2}&=&\sum
_{j=m}^{\infty} \biggl[\frac{C_{j}(\psi)}{j!}
\biggr]%
^{2}\sigma_{T}^{2}(j,z)
\nonumber
\\[-8pt]
\\[-8pt]
\nonumber
&=&\sum_{j=m}^{\infty}\frac{C_{j}^{2}(\psi)}{%
j!}\int
_{0}^{T}\int_{0}^{T}R_{T}(t)R_{T}(s)B^{j}(t-s)
\nu(dt)\nu(ds).
\end{eqnarray}
We will prove the asymptotic normality of (\ref{W}) under condition
(A4)(i). The proof under condition (A4)(ii) is even simpler.

By conditions (A2) and (A4)(i), for $j\geq2,$ all the
convolutions $f^{\ast(j)}$ are bounded and continuous functions,
and by (B1),
%
\begin{eqnarray} \label{lim4}
\sigma_{T}^{2}(j,z) &=&j!\int_{0}^{T}
\int_{0}^{T}B^{j}(t-s)R_{T}(t)R_{T}(s)
\nu(dt)\nu(ds)\nonumber
\\
&=&\sum_{k,l=1}^{q} \biggl( j!\int
_{0}^{T}\int_{0}^{T}B^{j}(t-s)
\frac{w_{k}(t)%
}{W_{k,T}}\frac{w_{l}(s)}{W_{l,T}}\nu(dt)\nu(ds) \biggr) z_{k}z_{l}
\nonumber
\\[-8pt]
\\[-8pt]
\nonumber
&=&2\pi j!\int_{\Lambda}f^{\ast(j)}(\lambda) \Biggl( \sum
_{k,l=1}^{q}\mu_{T}^{k,l}(
\lambda)z_{k}z_{l} \Biggr) \,d\lambda
\nonumber
\\
&\mathop{\longrightarrow}\limits_{T\rightarrow\infty}&2\pi j!\int_{\Lambda
}f^{\ast(j)}(
\lambda)m_{z}(d\lambda)=\sigma^{2}(j,z),\qquad j\geq2,
\nonumber
\end{eqnarray}
 where $m_{z}(d\lambda)=\sum_{k,l=1}^{q}\mu^{k,l}(d\lambda
)z_{k}z_{l}.$

Under condition (A2), from Theorem~\ref{atresul}, we obtain for $
j=1, $
%
\begin{equation}
\lim_{T\rightarrow\infty}\sigma_{T}^{2}(1,z)=2\pi\int
_{\Lambda
}f(\lambda)m_{z}(d\lambda)=
\sigma^{2}(1,z). \label{lim5}
\end{equation}
Thus%
%
\begin{equation}
\lim_{T\rightarrow\infty}\mathit{E} \langle\zeta_{T},z
\rangle^{2}=\sum_{j=1}^{\infty}
\biggl[\frac{C_{j}(\psi)}{j!}%
\biggr]^{2}\sigma^{2}(j,z)=
\sigma^{2}(z). \label{lim6}
\end{equation}
In Lemma~\ref{lemprevth} below, we will consider the following
decomposition:
%
\begin{eqnarray}\label{lim7}
\tau_{T}&=& \langle\zeta_{T},z \rangle=\tau_{T}(d)+
\tau_{T}^{\prime}(d)
\nonumber
\\[-8pt]
\\[-8pt]
\nonumber
&= &\Biggl( \sum_{j=1}^{d}+
\sum_{j=d+1}^{\infty} \Biggr) \frac{%
C_{j}(\psi)}{j!}
\int_{0}^{T}R_{T}(t)H_{j}
\bigl(\xi(t)\bigr)\nu(dt).
\end{eqnarray}

%
\begin{lemma}
\label{lemprevth} Suppose that conditions \textup{(A1)--(A4)} and
\textup{(B1)--(B3)} hold. If for any $d\geq1,$ as $T\rightarrow\infty,$
$\tau_{T}(d)\Rightarrow\tau_{d}\sim N(0,\sigma_{d}^{2}(z)),$
where
%
\begin{equation}
\sigma_{d}^{2}(z)=\sum_{j=1}^{d}
\biggl[\frac{C_{j}(\psi)}{j!} \biggr]%
^{2}\sigma^{2}(j,z),
\label{lim9}
\end{equation}
then $\tau_{T}\Rightarrow\tau\sim N(0,\sigma^{2}(z)).$
\end{lemma}

\begin{pf} Note that $\mathit{E}[(\tau_{T}^{\prime}(d))^{2}]\rightarrow
0, d\rightarrow\infty,$ uniformly in $T.$ Really, by condition
(B3),%
%
\begin{equation}
\bigl\llvert R_{T}(t)\bigr\rrvert=\Biggl\llvert\sum
_{i=1}^{q}z_{i}w_{i}(t)W_{i,T}^{-1}
\Biggr\rrvert\leq T^{-1/2}\llVert z\rrVert\llVert\tilde{k}\rrVert,\qquad
\tilde{k}=(k_{1},\ldots,k_{q}). \label{lim10}
\end{equation}
Then, under (A4)(i), as $d\rightarrow\infty,$
\begin{eqnarray*}
\mathit{E}\bigl[\bigl(\tau_{T}^{\prime}(d)\bigr)^{2}
\bigr] &=&\sum_{j=d+1}^{\infty}\frac{%
C_{j}^{2}(\psi)}{j!}
\int_{0}^{T}%
\int_{0}^{T}B^{j}(t-s)R_{T}(t)R_{T}(s)
\nu(dt)\nu(ds)
\\
&\leq&T^{-1}\llVert z\rrVert^{2} \llVert\tilde{k}
\rrVert^{2}\int_{0}^{T}\int
_{0}^{T}B^{2}(t-s)\nu(dt)\nu(ds)\sum
_{j=d+1}^{\infty}\frac{%
C_{j}^{2}(\psi)}{j!}
\\
&\leq&\llVert z\rrVert^{2} \llVert\tilde{k}\rrVert^{2}\int
_{\mathbb{R}}B^{2}(t)\nu(dt)\sum
_{j=d+1}^{\infty}\frac{C_{j}^{2}(\psi)}{j!}%
=\beta(d)
\longrightarrow0,
\end{eqnarray*}
since by Parseval's identity,%
\[
\sum_{j=1}^{\infty}\frac{C_{j}^{2}(\psi)}{j!}=
\mathit{E}\psi^{2}\bigl(\xi(0)\bigr)<\infty.
\]
Thus, for any $\varepsilon>0,$ uniformly in $T,$%
\[
P \bigl\{ \bigl\llvert\tau_{T}^{\prime}(d)\bigr\rrvert>
\varepsilon\bigr\} \leq\frac{\beta(d)}{\varepsilon
^{2}}\longrightarrow0, \qquad d\rightarrow\infty
.
\]

For any $\varepsilon>0,$ and $d\geq1,$ we then obtain
%
\begin{equation}
\overline{\lim_{T\rightarrow\infty}}P \{ \tau_{T}\leq x \} \leq
\Phi_{d}(x+\varepsilon)+\frac{\beta
(d)}{\varepsilon^{2}}, \label{lim11}
\end{equation}
 where $\Phi_{d}$ is the distribution function of a
Gaussian r.v. with zero mean and variance $\sigma_{d}^{2}(z).$

Also, for any $\varepsilon>0,$ and $d\geq1,$ as
$T\rightarrow\infty,$ the following inequality holds:
%
\begin{equation}
\underline{\lim}P \{ \tau_{T}\leq x \} \geq\Phi_{d}(x-
\varepsilon)-\frac{\beta(d)}{\varepsilon^{2}}. \label{lim12}
\end{equation}

If $d\rightarrow\infty,$ we obtain, from equations (\ref{lim11})
and (\ref{lim12}) that, as $T\rightarrow\infty,$
\[
\Phi_{\infty}(x-\varepsilon)\leq\underline{\lim}P \{ \tau_{T}
\leq x \} \leq\overline{\lim}P \{ \tau_{T}\leq x \} \leq
\Phi_{\infty}(x+\varepsilon),
\]
where $\Phi_{\infty}$ is the distribution function of a Gaussian
r.v. with zero mean and the variance $\sigma^{2}(z)$ given by (%
\ref{lim6}). Thus, if $\varepsilon\rightarrow
0, \lim_{T\rightarrow\infty}P \{ \tau_{T}\leq x \}
=\Phi_{\infty}(x),$ \mbox{$x\in\mathbb{R}$.}
\end{pf}

Now we are in position to derive the proof of Theorem
\ref{thmain2b}. In such a proof, we will check condition (i) of
Proposition~\ref{pr2}, but the proof can also be developed from the
verification of condition (ii) in Proposition~\ref{pr2}, using
diagram formula. We place this proof into \hyperref[app]{Appendix}, due to its
methodological interest in relation to the approach it presents for
the analysis of nonregular diagrams, providing the classification
of their levels into recipients and donors.

\subsection*{Proof of Theorem \protect\ref{thmain2b}}
From Lemma~\ref{lemprevth}, it is sufficient to show the asymptotic
normality of the r.v.'s $\tau_{T}(d)$. Consider then the r.v.'s%
%
\begin{equation}
\pi_{T,d}(\xi)= \biggl( \int_{0}^{T}r_{T,1}(t)H_{1}
\bigl(\xi(t)\bigr)\nu(dt),\ldots,\int_{0}^{T}r_{T,d}(t)H_{d}
\bigl(\xi(t)\bigr)\nu(dt) \biggr)^{\prime}, \label{lim13h}\hspace*{-35pt}
\end{equation}
 where
%
\begin{equation}
r_{T,j}(t)=\frac{R_{T}(t)}{\sigma(j,z)},\qquad j=1,\ldots,d. \label{lim13hh}
\end{equation}
The proof will follow from the application of Corollary~\ref{cor1},
after checking condition
(i) of Proposition~\ref{pr2} for the random vector $\pi_{T,d}(\xi)$
defined by (\ref{lim13h}) and~(\ref{lim13hh}). From Theorem 1 and equation
(\ref{lim4}),
%
\begin{eqnarray}\label{lim15}
\mathit{E} \biggl[ \int_{0}^{T}r_{T,j}(t)H_{j}
\bigl(\xi(t)\bigr)\nu(dt) \biggr]^{2}=%
\frac{\sigma_{T}^{2}(j,z)}{\sigma^{2}(j,z)}
\longrightarrow1,
\nonumber
\\[-8pt]
\\[-8pt]
\eqntext{  T\rightarrow\infty, j=1,\ldots,d.}
\end{eqnarray}

Now, $\pi_{T,d}(\xi)\Longrightarrow\pi_{d}\sim\mathcal{N} ( 0,%
\mathbb{I}_{d} ) ,$ $T\rightarrow\infty,$ if and only if
\[
\lim_{T\rightarrow\infty}\Vert f_{j,T}\otimes_{p}f_{j,T}
\Vert_{H^{\otimes2(j-p)}}=0
\]
 for $p=1,\ldots,j-1,$ $2\leq j\leq d,$ where
\[
f_{j,T}(s_{1},\ldots s_{j})=\int
_{0}^{T}R_{T}(t)\prod
_{i=1}^{j}\mathbf{l}%
_{(-\infty,t]}(s_{i})\,dt.
\]

We first check the convergence to zero of contractions in the
continuous time case. The $p$th contraction is computed by applying
formula (\ref{eqcontp})
with $k=j$ as follows:
%
\begin{eqnarray}\label{pcont}
&& f_{j,T}\otimes_{p} f_{j,T} (x_{1},
\ldots,x_{2j-2p})\nonumber\\[-2pt]
&&\qquad=\int_{0}^{T}\int
_{0}^{T}R_{T}(t)R_{T}(s)B(t-s)
\times\mathop{\cdots}_{p} \times B(t-s)\nonumber\\[-2pt]
&&\hspace*{68pt}{}\times \prod_{i=p+1}^{j}
\mathbf{l}_{(-\infty
,t]}(x_{i}) \prod_{l=p+1}^{j}
\mathbf{l}_{(-\infty,s]}(x_{l})\,ds\,dt
\\[-2pt]
&&\qquad=\int_{0}^{T}\int_{0}^{T}R_{T}(t)R_{T}(s)B^{p}(t-s)
\nonumber\\[-2pt]
&&\hspace*{67pt}{}\times\prod_{i=1}^{j-p} \mathbf{l}_{(-\infty,t]}(x_{i})\prod_{l=j-p+1}^{2j-2p}
\mathbf{l}%
_{(-\infty,s]}(x_{l})\,ds\,dt.\nonumber
\end{eqnarray}

The norm of the $p$th contraction (\ref{pcont}) in the space
$H^{\otimes2(j-p)}$ is then given by
\begin{eqnarray*}
&&\|f_{j,T}\otimes_{p}f_{j,T} \|^{2}_{H^{\otimes2(j-p)}}
\\[-2pt]
&&\qquad=\int_{0}^{T}\int_{0}^{T}
\int_{0}^{T}%
\int_{0}^{T}R_{T}(t_{1})R_{T}(s_{1})R_{T}(t_{2})R_{T}(s_{2})
\\[-2pt]
& & \hspace*{80pt}\qquad{} \times
B^{j-p}(t_{1}-t_{2})B^{j-p}(s_{1}-s_{2})B^{p}(t_{1}-s_{1})
\\[-2pt]
& &\hspace*{80pt}\qquad{} \times B^{p}(t_{2}-s_{2})\,ds_{1}\,ds_{2}\,dt_{1}\,dt_{2}.
\end{eqnarray*}
By condition (B3), for $2\leq j \leq d$ and
$p=1,\ldots,j-1,$
%
\begin{eqnarray}\label{ncg}
&& \|f_{j,T}\otimes_{p}f_{j,T}
\|^{2}_{H^{\otimes2(j-p)}}\nonumber\\
&&\qquad\leq \frac{%
\llVert z\rrVert^{4}\llVert \tilde{k}\rrVert^{4}}{T^{2}}%
\int_{0}^{T}
\int_{0}^{T}\int_{0}^{T}%
\int_{0}^{T}\bigl|B^{p}(t_{1}-s_{1})B^{p}(t_{2}-s_{2})\bigr|
\nonumber
\\
& &\hspace*{114pt}\qquad\quad {}\times\bigl|B^{j-p}(t_{1}-t_{2})B^{j-p}(s_{1}-s_{2})\bigr|
\,ds_{1}\,ds_{2}\,dt_{1} \,dt_{2}
\nonumber
\\
&&\qquad\leq \frac{\llVert z\rrVert^{4}\llVert
\tilde{k}\rrVert^{4}}{T^{2}}\int_{0}^{T}\int
_{0}^{T}\int_{0}^{T}%
\int_{0}^{T}\bigl|B(t_{1}-s_{1})B(t_{2}-s_{2})\bigr|
\nonumber
\\
& &\hspace*{114pt}\qquad\quad{} \times\bigl|B(t_{1}-t_{2})B(s_{1}-s_{2})\bigr|\,dt_{1}\,dt_{2}\,ds_{1}\,ds_{2}
\\
&&\qquad\leq\frac{\llVert z\rrVert^{4}\llVert \tilde{k}\rrVert^{4}%
}{T^{2}}\int_{0}^{T}\int
_{0}^{T}\int_{0}^{T}
\frac{1}{2}\int_{0}^{T} \bigl[%
B^{2}(t_{1}-s_{1})+B^{2}(t_{1}-t_{2})
\bigr]\,dt_{1}
\nonumber
\\
& &\hspace*{122pt}\qquad\quad{}  \times\bigl|B(t_{2}-s_{2})B(s_{1}-s_{2})\bigr|\,ds_{1}\,ds_{2}
\,dt_{2}
\nonumber
\\
&&\qquad\leq4\llVert z\rrVert^{4}\llVert\tilde{k}\rrVert^{4}%
\biggl[\int_{0}^{\infty}B^{2}(t_{1})\,dt_{1}
\biggr] \biggl[%
\int_{0}^{T}\bigl|B(t_{2})\bigr|\,dt_{2}
\biggr]
\nonumber
\\
& &\qquad\quad {} \times T^{-2} \int_{0}^{T}
\int_{0}^{T}\bigl|B(s_{1}-s_{2})\bigr|\,ds_{1}\,ds_{2}.\nonumber
\end{eqnarray}

In (\ref{ncg}), as $T\rightarrow\infty,$
%
\begin{eqnarray}\label{ordeint}
\int_{0}^{T}\bigl|B(t_{2})\bigr|\,dt_{2}
&=&\mathcal{O}\bigl(T^{1-\alpha}\bigr)
\nonumber
\\[-8pt]
\\[-8pt]
\nonumber
T^{-2}\int_{0}^{T}\int
_{0}^{T}\bigl|B(s_{1}-s_{2})\bigr|\,ds_{1}\,ds_{2}
&=&\mathcal{O}%
\bigl(T^{-\alpha}\bigr).
\end{eqnarray}
From condition (A4), in the case considered of Hermite
rank $m=1,$ we have $\alpha>1/2,$ and therefore, from (\ref{ncg})
and (\ref{ordeint}), we obtain for $j\geq2,$ $p=1,\ldots,j-1,$
%
\begin{equation}
\lim_{T\rightarrow\infty}\Vert f_{j,T}\otimes_{p}f_{j,T}
\Vert_{H^{\otimes2(j-p)}}=0. \label{concero}
\end{equation}

The proof in the discrete time case can be similarly derived in
terms of definition~(\ref{iso2}) of isonormal process $X,$
considering the counting measure $\nu(\cdot).$ Specifically, from
(\ref{dti}), for $T>0$ and $2\leq j\leq d$ we consider the
sequence of kernels
\[
f_{j,T}(m_{1},\ldots, m_{j})=\sum
_{l=1}^{T}R_{T}(l)\prod
_{i=1}^{j}\delta_{m_{i},l},\qquad  m_{1},
\ldots,m_{j}\in\mathbb{Z}.
\]

For $p=1,\ldots,j-1,$ the $p$th self-contraction of this kernel is given by
\begin{eqnarray*}
&& f_{j,T}\otimes_{p}
f_{j,T}(m_{1},\ldots,m_{2j-2p})\\
&&\qquad=\sum_{q=1}^{T}\sum
_{l=1}^{T}R_{T}(q)R_{T}(l)B^{p}(q-l)
\prod_{i=1}^{j-p} \delta_{m_{i},q} \prod_{i=j-p+1}^{2j-2p}
\delta_{m_{i},l}.
\end{eqnarray*}

Therefore, since
\begin{eqnarray*}
&&\|f_{j,T}\otimes_{p}f_{j,T} \|^{2}_{H^{\otimes2(j-p)}}\\
&&\qquad=
\sum_{q=1}^{T}\sum
_{k=1}^{T}\sum_{l=1}^{T}
\sum_{i=1}^{T} R_{T}(q)R_{T}(l)R_{T}(k)R_{T}(i)
\\
& &\hspace*{74pt}\qquad {}\times B^{j-p}(q-k)B^{j-p}(l-i)B^{p}(q-l)B^{p}(k-i),
\end{eqnarray*}
 in a similar way to the continuous time case, we obtain
%
\begin{eqnarray}\label{ncgb}
&&\|f_{j,T}\otimes_{p}f_{j,T} \|^{2}_{H^{\otimes2(j-p)}}
\nonumber\\
&&\qquad\leq \frac{\llVert z\rrVert^{4}\llVert \tilde{k}\rrVert^{4}}{T^{2}}
\sum_{q=1}^{T}
\sum_{k=1}^{T}\sum
_{l=1}^{T}\sum_{i=1}^{T}\bigl|B^{p}(q-l)B^{p}(k-i)\bigr|
\nonumber
\\
& &\hspace*{108pt}\qquad\quad{} \times\bigl|B^{j-p}(q-k)B^{j-p}(l-i)\bigr|
\nonumber
\\[-8pt]
\\[-8pt]
\nonumber
&&\qquad\leq \frac{\llVert z\rrVert^{4}\llVert
\tilde{k}\rrVert^{4}}{T^{2}}\sum_{q=1}^{T}
\sum_{k=1}^{T}\sum
_{l=1}^{T}\sum_{i=1}^{T}\bigl|B(q-l)B(k-i)\bigr|
\bigl|B(q-k)B(l-i)\bigr|\\
&&\qquad\leq4\llVert z\rrVert^{4}\llVert\tilde{k}
\rrVert^{4} \Biggl[\sum_{q=1}^{\infty}B^{2}(q)
\Biggr] \Biggl[%
\sum_{k=1}^{T}\bigl|B(k)\bigr|
\Biggr]
\nonumber\\
&&\qquad\quad{}\times  T^{-2} \sum_{l=1}^{T}
\sum_{i=1}^{T}\bigl|B(l-i)\bigr|.\nonumber
\end{eqnarray}

Thus, as $T\rightarrow\infty,$
%
\begin{eqnarray}\label{ordeintb}
\sum_{k=1}^{T}\bigl|B(k)\bigr| &=& \mathcal{O}
\bigl(T^{1-\alpha}\bigr),
\nonumber
\\[-8pt]
\\[-8pt]
\nonumber
T^{-2} \sum_{l=1}^{T}\sum
_{i=1}^{T}\bigl|B(l-i)\bigr| &=&\mathcal{O}
\bigl(T^{-\alpha}\bigr).
\end{eqnarray}

Again, from condition (A4), and equations
(\ref{ncgb}) and (\ref{ordeintb}), Proposition~\ref{pr2}(i) holds, and the
convergence to the Gaussian distribution follows.

\begin{example*}[(Continuation)]
Consider now model
(\ref{Reg1}) with nonlinear regression function
%
\begin{equation}
g(t,\theta)=\sum_{k=1}^{N} (
A_{k}\cos\varphi_{k}t+B_{k}\sin
\varphi_{k}t ) , \label{reg3}
\end{equation}
where $\theta=(A_{1},B_{1},\varphi_{1},\ldots,A_{N},B_{N},\varphi
_{N}), C_{k}^{2}=A_{k}^{2}+B_{k}^{2}>0, k=1,\ldots,N, 0<\varphi
_{1}<\cdots<\varphi_{N}<\infty.$ In this case, $q=3N,$ function $%
g(t,\theta)$ then has a block-diagonal measure $\mu(d\lambda)$
(see, e.g., Ivanov~\cite{r25}) with blocks
\[
\pmatrix{
\kappa_{k} & \mathrm{i}\rho_{k} &
\bar{\beta}_{k}
\vspace*{2pt}\cr
-\mathrm{i}\rho_{k}, & \kappa_{k} & \bar{
\gamma}_{k},
\vspace*{2pt}\cr
\beta_{k} & \gamma_{k}, & \kappa_{k},%
},\qquad  k=1,\ldots,N,
\]
where%
\[
\beta_{k}=\frac{\sqrt{3}}{2C_{k}} ( B_{k}\kappa_{k}+
\mathrm{i}%
A_{k}\rho_{k} ) ,\qquad
\gamma_{k}=\frac{\sqrt{3}}{2C_{k}} ( -A_{k}\kappa_{k}+
\mathrm{i}B_{k}\rho_{k} ).
\]
Here, the measure $\kappa_{k}=\kappa_{k}(d\lambda)$ and the
signed measure $\rho_{k}=\rho_{k}(d\lambda)$ are located at the
points $\pm\varphi_{k},$ and $\kappa_{k} ( \{ \pm
\varphi_{k} \} ) =\frac{1}{2}, \rho_{k} (
\{ \pm\varphi_{k} \}
) =\pm\frac{1}{2}.$ We then have%
\begin{eqnarray*}
g_{3k-2}(t,\theta)&=&\frac{\partial}{\partial A_{k}}g(t,\theta)=\cos
\varphi_{k}t,\qquad  g_{3k-1}(t,\theta)=\frac{\partial}{\partial B_{k}}%
g(t,\theta)=\sin\varphi_{k}t,
\\
g_{3k}(t,\theta)&=&\frac{\partial}{\partial\varphi_{k}}g(t,\theta)=-A_{k}t
\sin\varphi_{k}t+B_{k}t\cos\varphi_{k}t,\qquad  k=1,
\ldots,N.
\end{eqnarray*}
It is easy to see that if the s.d. $f$ satisfies (I), and
$\varkappa_{j}\neq\varphi_{k}, j=0,1,\ldots,r,$ $k=1,\ldots
,N,$ one can find a neighborhood $V_{j}(c_{0})$ of the point
$\varkappa_{j},$ for $j=0,1,\ldots,r,$ which does not contain the
points $\varphi_{k},$ $k=1,\ldots,N.$ Thus, for $T>T_{0},$ the
following condition holds:
\[
W_{iT}^{-1}\max_{\lambda\in V_{j}(c_{0})}\bigl\llvert
w_{T}^{i}(\lambda)\bigr\rrvert\leq k_{ij}T^{-1/2},\qquad
j\in J; i=3k-2,3k-1,3k; k=1,\ldots,N.
\]

In relation to the considered function $\mathbf{w}(t)=\nabla
g(t,\theta),$ the measure $\mu_{T}^{jl}(d\lambda)=\mu
_{T}^{jl}(d\lambda,\theta)$
approximates, in the weak sense, the spectral measure $\bolds{\mu}
(d\lambda)=\{\mu^{jl}(d\lambda)\}_{j,l=1}^{q}$ of the nonlinear
regression function $g(t,\theta)$ [see (\ref{Reg1})], where%
\begin{eqnarray*}
\mu_{T}^{jl}(d\lambda,\theta)&=&\frac{g_{T}^{j}(\lambda,\theta
)\overline{%
g_{T}^{l}(\lambda,\theta)}\,d\lambda}{\sqrt{\int_{\Lambda
}\llvert g_{T}^{j}(\lambda,\theta)\rrvert^{2}\,d\lambda
\int_{\Lambda}\llvert g_{T}^{l}(\lambda,\theta
)\rrvert^{2}\,d\lambda}},\qquad  j,l=1,
\ldots,q,
\\
g_{T}^{j}(\lambda,\theta)&=&\int_{0}^{T}e^{\mathrm{i}t\lambda
}g_{j}(t,
\theta)\nu(dt),\qquad j=1,\ldots,q
\end{eqnarray*}
 and $g_{j}(t,\theta)$ defines the $jth$ component of $\mathbf
{w(}%
t)=\nabla g(t,\theta),$ for $j=1,\ldots,q.$

If the s.d. $f(\lambda)$ satisfies condition (II), then $f$ is $
\bolds{\mu}$-admissible, and the block-diagonal matrix $\int
_{\Lambda}f(\lambda)%
\bolds{\mu}(d\lambda)$ consists of the blocks%
\[
f(\varphi_{k}) \pmatrix{ 1 & 0 & \displaystyle\frac{\sqrt{3}}{2}
\frac{B_{k}}{C_{k}}
\vspace*{2pt}\cr
0 & 1 & \displaystyle-\frac{\sqrt{3}}{2}\frac{A_{k}}{C_{k}}
\vspace*{2pt}\cr
\displaystyle\frac{\sqrt{3}}{2}\frac{B_{k}}{C_{k}} & \displaystyle-\frac{\sqrt{3}}{2}\frac{A_{k}}{C_{k}%
} &
1%
},\qquad  k=1,\ldots,N.
\]

It is easy to see that, for the function $g(t,\theta)$ given by (\ref
{reg3}%
), the matrix $\Xi$ is a block-diagonal with blocks of the form%
\begin{eqnarray}
\Xi_{k}=2\pi\sum_{j=m}^{\infty}
\frac{C_{j}^{2}(\psi
)}{j!}f^{\ast(j)}(\varphi_{k}) \pmatrix{ 1 & 0 & \displaystyle B_{k}\frac{\sqrt{3}}{2C_{k}}
\vspace*{2pt}\cr
0 & 1 & \displaystyle -A_{k}\frac{\sqrt{3}}{2C_{k}}
\vspace*{2pt}\cr
\displaystyle B_{k}\frac{\sqrt{3}}{2C_{k}} & \displaystyle -A_{k}\frac{\sqrt{3}}{2C_{k}} &
1} ,
\nonumber\\
 \eqntext{ k=1,\ldots,N.}
\end{eqnarray}
\end{example*}

%
\section{Final comments}\label{sec6}

This paper addresses the problem of Gaussian limit theory of
weighted functionals of nonlinear transformations of Gaussian
stationary random processes $\xi$ having multiple singularities in
their spectra. The general case where the Fourier transform of the
weight function also displays multiple singularities in the limit,
which do not coincide with the singularities of the spectral density
of $\xi,$ is also covered here. This subject has several
applications in asymptotic statistical inference. We are especially
motivated by its application in the limit theory of nonlinear
regression problems with regression function and errors having
multiple singularities in their spectra. This actually constitutes
an active research area, due to the existence of several open
problems and applications. Note that, although here we have
considered the parameter range
\[
\alpha=\min_{j=0,1,\ldots,r}\alpha_{j}>1/2,
\]
 which, in particular, allows us to consider long-range
dependence models. Our conjecture is that the Gaussian limit results
hold for $\alpha_{j}\in(0,1),$ $j=0,1,\ldots,r.$ The proof of this
conjecture will lead to a general scenario where most of the limit
results derived for random fields with singular spectra (see Taqqu
\cite{r46,r47}; Dobrushin and Major~\cite{r15}; Nualart and Peccati
\cite{r36}; and the references therein) can be obtained as
particular cases.

\begin{appendix}

\section*{\texorpdfstring{Appendix: Proof of Theorem \lowercase{\protect\ref{thmain2b}} based on diagram formula}
{Appendix: Proof of Theorem 5.1 based on diagram formula}}\label{app}

As before, we will prove this result for Hermite rank $m=1.$ To show
the asymptotic normality of the r.v.'s $\tau_{T}(d),$ consider the
r.v.'s $\pi_{T,d}(\xi)$ and $r_{T,j}(t),$ $j=1,\ldots, d,$ defined
by (\ref{lim13h}) and (\ref{lim13hh}).

We will check condition (ii) of Proposition~\ref{pr2}. Then, from
Corollary~\ref{cor1}, $\pi_{T,d}(\xi
)\Rightarrow\pi_{d}\sim N(0,\mathbb{I}_{d})$, that is, $\tau_{T}(d)%
\Rightarrow\tau_{d}\sim N(0,\sigma_{d}^{2}(z)),$ as $T\rightarrow
\infty.$

We apply diagram technique for proving condition (ii) of Proposition
\ref{pr2}. Let us first introduce some definitions.\vadjust{\goodbreak}

A graph $\Gamma=$ $\Gamma(l_{1},\ldots,l_{p})$ with $l_{1}+\cdots
+l_{p}$ vertices is called a diagram of order $(l_{1},\ldots
,l_{p})$ if:
\begin{longlist}[(a)]
\item[(a)] the set of vertices $V$ of the graph $\Gamma$ is of the
form $%
V=\bigcup_{j=1}^{p}W_{j}$, where $W_{j}=\{(j,l)\dvtx1\leq l\leq l_{j}\}$
is the $j$th level of the graph $\Gamma$, $1\leq j\leq p$ (if $l_{j}=0,$ assume $
W_{j}=\varnothing$);

\item[(b)] each vertex is of degree 1;

\item[(c)] if $((j_{1},l_{1}), (j_{2},l_{2}))\in\Gamma, $ then
$j_{1}\neq j_{2}$, that is, the edges of the graph $\Gamma$ may
connect only different levels.
\end{longlist}

Let $L=L(l_{1},\ldots,l_{p})$ be a set of diagrams $\Gamma$ of order $
(l_{1},\ldots,l_{p})$. Denote by $Z_{\Gamma}$ the set of edges of a
graph $%
\Gamma\in L$. For the edge $\varpi=((j_{1},l_{1}), (j_{2},l_{2}))\in
Z_{\Gamma}$, $ j_{1}<j_{2}$, we set $d_{1}(\varpi)=j_{1} $,$%
d_{2}(\varpi)=j_{2}$. We call a diagram $\Gamma$ regular if its
levels can be split into pairs in such a manner that no edge
connects the levels belonging to different pairs. We denote by
$L^{\ast}$ the set of regular
diagrams $L^{\ast}\subseteq L(l_{1},\ldots,l_{p})$. If $p$ is odd,
then $%
L^{\ast}=\varnothing$.

The following lemma provides the diagram formula; see Taqqu
\cite{r47}, Lem\-ma~3.2 or Doukhan, Oppenheim and Taqqu~\cite{r17},
page~74, or Peccati and Taqqu~\cite{r39}.

%
\begin{lemma}
Let $(\xi_{1},\ldots,\xi_{p}), p\geq2$, be a Gaussian vector with $%
\mathit{E}\xi_{j}=0,\break \mathit{E}\xi_{j}^{2}=1, \mathit{E}\xi_{i}\xi
_{j}=B(i,j), i,j=1,\ldots,p,$ and let $H_{l_{1}}(u),\ldots
,H_{l_{p}}(u)$ be the Hermite polynomials. Then
%
\begin{equation}
\mathit{E} \Biggl\{ \prod_{j=1}^{p}H_{l_{j}}(
\xi_{j}) \Biggr\} =\sum_{\Gamma\in L}\prod
_{\varpi\in Z_{\Gamma}}B\bigl(d_{1}(\varpi
),d_{2}(\varpi)\bigr). \label{lim17}
\end{equation}
\end{lemma}

From (\ref{lim17}), we obtain, for $p=4,$ $%
l_{1}=l_{2}=l_{3}=l_{4}=j,$ $\Gamma=\Gamma(j,j,j,j)$ and $(\xi_{1},\xi
_{2},\xi_{3},\xi_{4})=(\xi(t_{1}),\xi(t_{2}),\xi
(t_{3}),\xi(t_{4})),$
%
\begin{eqnarray}\label{lim19}
&&\mathit{E}\pi_{T,j}^{4}(\xi)\nonumber\\
&&\qquad=\int_{0}^{T}
\int_{0}^{T}\int_{0}^{T}
\int_{0}^{T}%
\prod
_{i=1}^{4}r_{T,j}(t_{i})\\
&&\hspace*{102pt}{}\times
\mathit{E} \Biggl[\prod_{i=1}^{4}H_{j}
\bigl(\xi(t_{i})\bigr) \Biggr]\nu(dt_{1})\nu
(dt_{2})\nu(dt_{3})\nu(dt_{4}).\nonumber
\end{eqnarray}

We then have
\begin{eqnarray}\label{lim20}\qquad
\mathit{E}\pi_{T,1}^{4}(\xi)&=&\frac{3}{\sigma^{4}(1,z)}
\biggl[%
\int_{0}^{T}\int
_{0}^{T}B(t_{1}-t_{2})R_{T}(t_{1})R_{T}(t_{2})
\nu(dt_{1})\nu(dt_{2}) \biggr]^{2}
\nonumber
\\[-8pt]
\\[-8pt]
\nonumber
&=& 3\frac{\sigma_{T}^{4}(1,z)}{\sigma^{4}(1,z)}\longrightarrow3,\qquad
T\rightarrow\infty.
\end{eqnarray}

For $j\geq2,$ the sum in (\ref{lim17}) is split into two sums
corresponding to regular and nonregular diagrams,
\[
\sum_{\Gamma\in L}\cdots=\sum_{\Gamma\in L^{\ast}}\cdots+
\sum_{\Gamma
\in L\setminus L^{\ast}}\cdots,
\]
and the right-hand side of (\ref{lim19}) is split into
two these parts, as well.

\textit{Analysis of the regular diagrams}:

We have
%
\begin{equation}
\sum^{\ast} ( T ) =\sum_{\Gamma\in L^{\ast}}F_{\Gamma
}(T),
\label{lim21}
\end{equation}
 where
%
\begin{eqnarray}\label{fgamma}\qquad
&&F_{\Gamma
}(T)\nonumber\\
&&\qquad=\int_{0}^{T}\int
_{0}^{T}\int_{0}^{T}
\int_{0}^{T}%
\prod
_{i=1}^{4}r_{T,j}(t_{i})\\
&&\hspace*{103pt}{}\times \prod
_{\varpi\in Z_{\Gamma
}}B(t_{d_{1}(\varpi)}-t_{d_{2}(\varpi)})\nu
(dt_{1})\nu(dt_{2})\nu(dt_{3})\nu
(dt_{4}).\nonumber
\end{eqnarray}

Each regular diagram $\Gamma\in L^{\ast}$ consists of $ 4$ levels
of cardinality $j.$ There are only $3$ subdivisions of the $4$
levels into pairs, and in each pair the vertices can be connected by
$j!$ ways. Thus, there is only
\[
\bigl\llvert L^{\ast}\bigr\rrvert=3(j!)^{2}
\]
regular diagrams, and, in this case, sum (\ref{lim21}) is
subdivided into product of pairs of integrals
%
\begin{eqnarray}\label{lim22}\qquad
 \sum^{\ast} ( T ) &=&\frac{3(j!)^{2}}{\sigma^{4}(j,z)} \biggl(
\int_{0}^{T}\int_{0}^{T}B^{j}(t_{1}-t_{2})R_{T}(t_{1})R_{T}(t_{2})
\nu(dt_{1})\nu(dt_{2}) \biggr)^{2}
\nonumber
\\[-8pt]
\\[-8pt]
\nonumber
&  =& 3\frac{\sigma_{T}^{4}(j,z)}{\sigma
^{4}(j,z)}\longrightarrow3,\qquad  T\rightarrow\infty.
\end{eqnarray}

\textit{Analysis of the nonregular diagrams:}

First, we consider
%
\begin{equation}
\sum( T ) =\sum_{\Gamma\in L\setminus L^{\ast
}}F_{\Gamma}(T),
\label{eq111}
\end{equation}
 where $F_{\Gamma}$ is defined as in (\ref{fgamma}). We
now prove that $\lim_{T\rightarrow\infty}\sum( T ) =0.$
Then,
the assertion of the theorem will follow from (\ref{lim20}) and (%
\ref{lim22}).

From (\ref{lim10}),
\begin{eqnarray}\label{lim25}
\qquad\bigl\llvert F_{\Gamma}(T)\bigr\rrvert
&\leq&\frac{\llVert
z\rrVert^{4}\llVert \tilde{k}\rrVert^{4}}{\sigma^{4}(j,z)}%
T^{-2}
\nonumber\\
& &{} \times\int_{0}^{T}\int
_{0}^{T}\int_{0}^{T}
\int_{0}^{T}\prod_{i=1}^{4}
\prod_{%
\varpi\in Z_{\Gamma},d_{1}(\varpi)=i}\bigl\llvert
B(t_{i}-t_{d_{2}(\varpi)})
\bigr\rrvert\\
&&\hspace*{158pt}{}\times\nu(dt_{1})\nu(dt_{2})\nu(dt_{3})
\nu(dt_{4}).\nonumber
\end{eqnarray}

Let $q_{\Gamma}(i)$  be the number of edges $\varpi\in Z_{\Gamma
},$ such that $d_{1}(\varpi)=i.$ Then, for $q_{\Gamma}(i)\geq1,$
%
\begin{eqnarray}\label{lim26}
& & \int_{0}^{T}\prod
_{\varpi\in Z_{\Gamma},d_{1}(\varpi
)=i}\bigl\llvert B(t_{i}-t_{d_{2}(\varpi)})\bigr
\rrvert\nu(dt_{i})
\nonumber
\\
& & \qquad\leq\frac{1}{q_{\Gamma
}(i)}\sum_{\varpi\in Z_{\Gamma},d_{1}(\varpi
)=i}
\int_{0}^{T}\bigl\llvert B(t_{i}-t_{d_{2}(\varpi
)})
\bigr\rrvert^{q_{\Gamma}(i)}\nu(dt_{i})
\\
& &\qquad \leq2\int_{0}^{T}\bigl\llvert
B(t_{i})\bigr\rrvert^{q_{\Gamma}(i)}\nu(dt_{i}).\nonumber
\end{eqnarray}

If $q_{\Gamma}(i)$ $=0,$ the integrals regarded to these variables
($t_{4},$ and possibly $t_{3}),$ in the left-hand side of
(\ref{lim26}),
give a contribution in the form of a multiplier of $T$ in the estimate
(\ref%
{lim25}).

%
\begin{definition}
The level $i$ of a nonregular diagram $\Gamma\in L\setminus
L^{\ast}$ is said to be a donor, if $q_{\Gamma}(i)$ $\geq1,$
and a strong donor, if $q_{\Gamma}(i)=j.$ The level $i$ of a
nonregular diagram $\Gamma\in L\setminus L^{\ast}$ is said to be
a recipient, if it is not donor, that is $q_{\Gamma}(i)=0.$
\end{definition}

Let $\rho_{sd}$ be a number of strongly donor levels, and
$\rho_{r} $ be a number of recipient levels. Obviously, level
$1$ is a strong donor, while level $4$ is a recipient.
If $\rho_{\mathrm{sd}}=1,$ then $\rho_{r}=1,$ while if $\rho_{\mathrm{sd}}=2,$ then $\rho_{r}=2.$

Formulas (\ref{lim25}) and (\ref{lim26}) then imply
%
\begin{equation}
\bigl\llvert F_{\Gamma}(T)\bigr\rrvert\leq2^{4-\rho
_{r}}
\frac{\llVert
z\rrVert^{4}\llVert \tilde{k}\rrVert^{4}}{\sigma^{4}(j,z)}%
T^{-2}\prod_{i=1}^{4}
\int_{0}^{T}\bigl\llvert B(t)\bigr
\rrvert^{q_{\Gamma}(i)}\nu(dt). \label{lim27}
\end{equation}

Since $j\geq2,$ and $\alpha>1/2,$ for a strong donor level $i$
with $q_{\Gamma}(i)=j$,
%
\begin{equation}
\int_{0}^{T}\bigl\llvert B(t)\bigr
\rrvert^{j}\nu(dt)\leq\int_{0}^{\infty}
\bigl[B(t)\bigr]^{2}\nu(dt)<\infty. \label{lim28}
\end{equation}

Thus, for the recipient levels ($q_{\Gamma}(i)=0)$ and the strong
donor
levels ($q_{\Gamma}(i)=j$), we obtain%
%
\begin{equation}
\int_{0}^{T}\bigl\llvert B(t)\bigr
\rrvert^{q_{\Gamma}(i)}\nu(dt)\leq C_{0}T^{1-z(i)}, \label{lim29}
\end{equation}
 where
\[
z(i)=\frac{q_{\Gamma}(i)}{j}, \qquad C_{0}=\max\biggl( 1,\int
_{0}^{\infty}B^{2}(t)\nu(dt) \biggr).
\]

Let now $0<q_{\Gamma}(i)<j;$ that is, level $i$ is a donor,
but not strong donor, and then
\begin{eqnarray} \label{lim30}
\int_{0}^{T}\bigl\llvert B(t)\bigr
\rrvert^{q_{\Gamma}(i)}\nu(dt)&=& \biggl[%
\int_{0}^{1}+
\int_{1}^{T} \biggr]\bigl\llvert B(t)\bigr
\rrvert^{q_{\Gamma}(i)}\nu(dt)
\nonumber
\\
&\leq& 1+\frac{T^{1-\alpha q_{\Gamma}(i)}-1}{1-\alpha q_{\Gamma
}(i)}\\
&=&\frac{%
\alpha q_{\Gamma}(i)}{\alpha q_{\Gamma}(i)-1}+\frac{T^{1-\alpha
q_{\Gamma}(i)}}{1-\alpha q_{\Gamma}(i)}
=o \bigl(
T^{1-z(i)} \bigr) ,\nonumber
\end{eqnarray}
since $\alpha q_{\Gamma}(i)=\alpha jz(i),$ and $\alpha j>1.$ We
will show that
\[
\mu=2-\sum_{i=1}^{4}z(i)=0.
\]

Indeed,
\[
\sum_{i=1}^{4}z(i)=1+\frac{q_{\Gamma}(2)+q_{\Gamma}(3)}{j},
\]
and $q_{\Gamma}(2)+q_{\Gamma}(3)=j,$ since $\llvert
Z_{\Gamma}\rrvert =2j.$

Formulas (\ref{lim29}), (\ref{lim30}) and (\ref{lim26}) together with
(\ref%
{lim27}) then imply that
%
\begin{equation}
\bigl\llvert F_{\Gamma}(T)\bigr\rrvert=O(1),\qquad T\rightarrow\infty,
\label{lim31}
\end{equation}
when $\rho_{\mathrm{sd}}=$ $\rho_{r}=2,$ and
%
\begin{equation}
\bigl\llvert F_{\Gamma}(T)\bigr\rrvert=o(1), \qquad T\rightarrow\infty,
\label{lim32}
\end{equation}
when $0<q_{\Gamma}(i)<j,$ for $i=2,3$ ($\rho_{\mathrm{sd}}=$ $\rho_{r}=1$).

The estimate (\ref{lim31}) is not exact. Thus, let us consider again
the case of nonregular diagram $\Gamma,$ which has $2$ strong
donor levels, and the remaining $2$ levels are recipients. The
recipient level $3$ takes edges from the strong donor levels $1$
and~$2,$ while level $ 2$ does not supply level $3$ in
full. Let us permutate levels $2$ and $3,$ and denote this
permutation by $\pi,$ that is, $\pi(2)=3,\break\pi(3)=2,$ and, from the
level $\pi(3)$ to the level $\pi(2),$ there are less than $j$
edges. Moreover, from the level $\pi(3)$ there is no edges down,
except the edges which connect $\pi(3)$ with $\pi(2),$ since
level $\pi(3)$ took\vadjust{\goodbreak}
all edges from the top, that is,%
\[
q_{\pi\Gamma}\bigl(\pi(3)\bigr)=q_{\pi\Gamma}(2)<j,
\]
where $\pi\Gamma$ is a nonregular diagram, taken from $\Gamma$ by
permutating the levels $2$ and~$3.$ Note that this permutation does
not change the value of integral defining $F_{\Gamma}(T)$ in~(\ref{eq111}),
since it is equivalent to the renaming of the
variables $t_{2}$
and $t_{3}.$ From~(\ref{lim32}), we then obtain%
%
\begin{equation}
\bigl\llvert F_{\Gamma}(T)\bigr\rrvert=\bigl\llvert F_{\pi\Gamma
}(T)
\bigr\rrvert\rightarrow0, \qquad T\rightarrow\infty. \label{lim33}
\end{equation}

The assertion of this theorem then follows from equations~(\ref
{lim31})--(\ref{lim33}).
\end{appendix}

\section*{Acknowledgements}
The authors are grateful to Professor A. Yu. Pilipenko for helpful
discussions. The authors wish to thank the anonymous referees for
many useful comments and suggestions.

%


\printaddresses

\end{document}